\documentclass{article}
\usepackage{arxiv}
\usepackage{enumitem}
\usepackage[utf8]{inputenc}
\usepackage[T1]{fontenc}
\usepackage{amsmath,amssymb,amsthm,mathtools}
\usepackage{microtype}
\usepackage{algorithm}
\usepackage{algpseudocode}
\usepackage{booktabs,tabularx}
\usepackage[square,numbers]{natbib}
\usepackage{xcolor}
\usepackage{listings}
\lstdefinestyle{pythonpaper}{
    language=Python,
    basicstyle=\ttfamily\footnotesize,
    keywordstyle=\color{blue!60!black},
    commentstyle=\color{black!55},
    stringstyle=\color{black!70},
    numbers=none,
    frame=lines,
    rulecolor=\color{black!25},
    columns=fullflexible,
    keepspaces=true,
    showstringspaces=false,
    breaklines=true,
    captionpos=b
}
\usepackage{tikz}
\usetikzlibrary{arrows.meta,decorations.pathreplacing}
\usepackage{placeins}
\definecolor{linkcolor}{RGB}{178,129,0}
\usepackage[colorlinks=true,allcolors=linkcolor]{hyperref}
\usepackage[capitalize,nameinlink]{cleveref}

\theoremstyle{plain}
\newtheorem{theorem}{Theorem}[section]
\newtheorem{lemma}[theorem]{Lemma}

\newtheorem{corollary}{Corollary}[section]   % own counter
\theoremstyle{definition}
\newtheorem{assumption}[theorem]{Assumption}
\newtheorem{remark}{Remark}[section]         % bold head, own counter
\crefname{assumption}{assumption}{assumptions}
\Crefname{assumption}{Assumption}{Assumptions}

\newcommand{\R}{\mathbb{R}}

\newcommand{\half}{\frac{1}{2}}
\newcommand{\iprod}[2]{\left\langle #1, #2 \right\rangle}
\newcommand{\iprodi}[2]{\langle #1, #2 \rangle}

\newcommand\tr{\mathrm{tr}}

\newcommand{\order}[1]{\mathcal{O}\left(#1\right)}

\newcommand{\prt}[1]{\left(#1\right)}
\newcommand{\brk}[1]{\left[#1\right]}
\newcommand{\crk}[1]{\left\{#1\right\}}

\newcommand{\orderi}[1]{\mathcal{O}(#1)}
\newcommand{\prti}[1]{(#1)}
\newcommand{\brki}[1]{[#1]}
\newcommand{\crki}[1]{\{#1\}}

\newcommand{\normf}[1]{\left\Vert #1 \right\Vert_{\mathrm{F}}}
\newcommand{\normif}[1]{\Vert #1 \Vert_{\mathrm{F}}}

\newcommand{\norm}[1]{\left\Vert #1 \right\Vert}
\newcommand{\normi}[1]{\Vert #1 \Vert}

\newcommand{\x}{\mathbf{x}}

\newcommand{\1}{\mathbf{1}}

\newcommand{\T}{\intercal}

\newcommand{\s}{\mathbf{s}}

\DeclareMathOperator{\Pol}{P}
\DeclareMathOperator{\sgn}{sign}

\definecolor{cuhkpl}{RGB}{152,24,147}
\definecolor{cuhkp}{RGB}{255,1,255}

\title{Muon Attains $\mathcal{O}(\varepsilon^{-3})$ Iteration Complexity under Shuffling}

\author{
Kun Huang \\
KTH Royal Institute of Technology \\
Department of Decision and Control Systems (DCS) \\
Digital Futures\\
Stockholm, Sweden \\
\texttt{kunhuang@kth.se}
\And
Shi Pu\\
The Chinese University of Hong Kong, Shenzhen \\
School of Data Science (SDS) \\
Shenzhen, Guangdong, China \\
\texttt{pushi@cuhk.edu.cn}
\And
Karl Henrik Johansson \\
KTH Royal Institute of Technology \\
Department of Decision and Control Systems (DCS) \\
Digital Futures\\
Stockholm, Sweden \\
\texttt{kallej@kth.se}
}

\begin{document}
\raggedbottom
\maketitle

\begin{abstract}
Muon is a stochastic first-order optimizer for matrix-valued parameters. It maintains an exponential moving average of stochastic gradients in a momentum buffer with coefficient $\beta$, and orthogonalizes the buffer to obtain the update direction. Muon has been reported to accelerate and stabilize the training of large language models. In practice, it is often used with shuffling (sampling without replacement), where each epoch visits the components of a finite-sum objective in a permuted order. However, its iteration complexity in this setting remains open.
We analyze Muon under shuffling, with the momentum buffer orthogonalized by the polar map. For smooth nonconvex finite-sum problems, we establish an $\mathcal{O}(\varepsilon^{-3})$ iteration complexity for finding an $\varepsilon$-stationary point. This guarantee holds for every sequence of permutations, covering random reshuffling, shuffle once, and the incremental gradient method.
This complexity matches that of stochastic gradient descent under shuffling in its dependence on $\varepsilon$, and improves on Muon's $\mathcal{O}(\varepsilon^{-4})$ complexity under sampling with replacement. 
Our guarantee requires a problem-dependent $\beta$ that approaches one as the number of components grows. In contrast, no fixed $\beta$ guarantees  convergence to stationarity for every number of components.
We show this by constructing a family of smooth convex problems on which, for any fixed $\beta\in[0,1)$, Muon remains bounded away from stationarity once the number of components is sufficiently large. This failure occurs for every sequence of permutations and persists even with small stepsizes.
These results suggest two practical tuning rules: if the training loss stops decreasing even after reducing the stepsize, increase $\beta$; when increasing $\beta$, reduce the stepsize accordingly.
Numerical experiments are consistent with our theoretical findings and support both recommendations.

\end{abstract}

\section{Introduction}\label{sec:intro}

Many machine learning models are trained by minimizing an empirical risk over matrix-valued parameters:
\begin{equation}\label{eq:P}
    \min_{W\in\R^{m\times n}} f(W)=\frac{1}{N}\sum_{i=1}^{N}f_i(W),
\end{equation}
where $W$ is a weight matrix, $N$ is the number of components, and each $f_i$ is the loss on a single data sample or a mini-batch from a fixed data partition.
Classical optimizers, including stochastic gradient descent with momentum (SGDM)~\cite{polyak1964some}, AdaGrad~\cite{duchi2011adaptive}, and Adam~\cite{kingma2014adam}, typically treat $W$ as a vector and thus do not exploit its matrix structure~\cite{bernstein2024old,gupta2018shampoo}. Muon (MomentUm Orthogonalized by Newton--Schulz)~\cite{jordan2024muon}, by contrast, is designed to leverage this structure: it maintains a momentum buffer $M$ and updates the parameters along an orthogonalization of $M$. Muon has been adopted to train large language models such as DeepSeek-V4~\cite{deepseekv4} and Kimi K3~\cite{kimik3}, whose developers report faster convergence and more stable training.
Like SGDM and Adam, Muon is typically used with stochastic gradients, where each update uses the gradient of a single component. Consequently, its convergence depends on how the components are sampled.

Most existing convergence analyses of Muon~\cite{li2025note,shen2026on,wang2026convergence} assume sampling with replacement, in which a component is drawn independently and uniformly at random at each update.
Under sampling without replacement (also known as shuffling), each component is visited exactly once per epoch, in the order specified by a permutation. 
Shuffling schemes differ in how the permutations are chosen: random reshuffling (RR)~\cite{li2023convergence,mishchenko2020random} draws a fresh uniformly random permutation at each epoch; shuffle once (SO)~\cite{safran2020good} draws one random permutation and reuses it across epochs; and the incremental gradient (IG) method~\cite{bertsekas2011incremental} uses a fixed, deterministic permutation. 

In practice, Muon is typically used with shuffling. For example, the official Moonlight training example combines Muon with a shuffled data loader~\cite{liu2025muonscalablellmtraining}.\footnote{In Moonshot AI's script \href{https://github.com/MoonshotAI/Moonlight/blob/c2ad5b20c605086526a179d36901bfc41b52b44b/examples/toy_train.py\#L242-L306}{\texttt{examples/toy\_train.py}}, the function \texttt{get\_model\_and\_dataloader} builds a \texttt{DataLoader} with \texttt{shuffle=True}, and \texttt{get\_optimizer} applies Muon to the hidden weight matrices.}
Motivated by this practice, we study Muon under shuffling. To cover RR, SO, and the IG method within a unified analysis, we allow an arbitrary sequence of permutations $\{\pi_k\}_{k\geq0}$. 
At step $j$ of epoch $k$, corresponding to iteration $t=kN+j$, Muon updates its momentum buffer and parameters as follows:
\begin{equation}\label{eq:intro-muon}
    M_t=\beta M_{t-1}+(1-\beta)\nabla f_{\pi_k(j+1)}(W_t),\quad
    W_{t+1}=W_t-\eta_t\Pol(M_t), \qquad k\geq0,\quad j\in\{0,\ldots,N-1\},
\end{equation}
where $\eta_t>0$ are the stepsizes and the momentum coefficient $\beta\in[0,1)$ is constant over iterations. The momentum buffer is carried across epochs. The orthogonalization is given by the polar map $\Pol(M)=U_sV_s^{\T}$, where $M=U_s\Sigma_s V_s^{\T}$ is a compact singular value decomposition (SVD), with $\Pol(\mathbf{0}) = \mathbf{0}$; see \eqref{eq:polar}. We analyze exact-polar Muon, whereas practical implementations approximate it using a few Newton--Schulz iterations~\cite{jordan2024muon}. 

\Cref{tab:comparison} summarizes known iteration complexities for smooth nonconvex problems. Under sampling with replacement, SGD~\cite{ghadimi2013stochastic}, SGDM~\cite{liu2020improved}, Adam~\cite{jin2024comprehensive}, and Muon~\cite{li2025note,shen2026on,wang2026convergence} admit an $\orderi{\varepsilon^{-4}}$ iteration complexity. Under shuffling, SGD and SGDM achieve an improved $\orderi{\varepsilon^{-3}}$ iteration complexity on $\varepsilon$ for finite-sum problems satisfying suitable growth conditions~\cite{mishchenko2020random,nguyen2020unified,qiu2024random}. Adam, in contrast, converges to a neighborhood of stationary points whose radius vanishes as its second-moment coefficient $\beta_2$ approaches one~\cite{zhang2026adam}. For Muon under RR, recent studies have focused on its implicit bias~\cite{li2026implicit} and generalization~\cite{kim2026geometry}, rather than its iteration complexity.

\begin{table}[htbp]
\centering
\caption{Iteration complexity for finding an $\varepsilon$-stationary point in smooth nonconvex
problems.
Stationarity is measured in the nuclear norm for Muon and in the Euclidean norm for the other methods.}
\label{tab:comparison}
\begin{tabular}{@{}lcc@{}}
\toprule
     & Sampling with replacement & Shuffling \\ \midrule
SGD  & $\order{\varepsilon^{-4}}$~\cite{ghadimi2013stochastic} & $\order{\varepsilon^{-3}}$~\cite{mishchenko2020random,nguyen2020unified} \\
SGDM & $\order{\varepsilon^{-4}}$~\cite{liu2020improved}       & $\order{\varepsilon^{-3}}$~\cite{qiu2024random} \\
Adam & $\order{\varepsilon^{-4}}$~\cite{jin2024comprehensive}  & neighborhood$^{\dagger}$~\cite{zhang2026adam} \\
Muon & $\order{\varepsilon^{-4}}$~\cite{li2025note,shen2026on,wang2026convergence} & $\boldsymbol{\order{\varepsilon^{-3}}}$ \textbf{(this work)}$^{\ddagger}$ \\
\bottomrule
\end{tabular}

\smallskip
\begin{minipage}{\linewidth}\footnotesize
$^{\dagger}$\citet{zhang2026adam} prove a bound of
$\orderi{\log T/\sqrt{T}+\delta(\beta_2)\sqrt{D_0}}$, where $D_0\geq0$ is defined in \eqref{eq:growth-intro}, $\beta_2$ is Adam's second-moment coefficient, and $\delta(\beta_2)\to0$ as $\beta_2\to1$ at an unspecified rate.\\
$^{\ddagger}$ This holds for every sequence of permutations, including RR, SO, and IG.
\end{minipage}
\end{table}

Typical shuffling analyses for SGD do not directly extend to Muon. They rely on the fact that, over each epoch, the component gradients average to the full gradient~\cite{li2023convergence,huang2023drr}. This argument requires the update to depend linearly on the gradients, whereas Muon applies the nonlinear polar map. Momentum may restore this averaging before the polar map: the momentum buffer averages roughly $1/(1-\beta)$ recent component gradients. Thus, it can approximate the full gradient only when $1/(1-\beta)$ is at least of order $N$. This suggests choosing $\beta$ as a function of $N$, whereas in practice $\beta$ is often set to a default value such as $0.95$~\cite{jordan2024muon}. These observations motivate the following question:
\begin{quote}
\centering
\emph{Can Muon under shuffling attain an $\orderi{\varepsilon^{-3}}$ iteration complexity, and how should the momentum coefficient be chosen to guarantee convergence?}
\end{quote}

We answer both parts: Muon under shuffling attains an $\orderi{\varepsilon^{-3}}$ iteration complexity with a suitable stepsize and a problem-dependent $\beta$ that approaches one as the number of components grows.
Our analysis relies on three standard assumptions: each component has a Lipschitz continuous gradient; $f$ is bounded below by $f^* = \inf_W f(W)$; and, for some constants $D_1,D_0\geq0$, the component gradients satisfy the following growth condition:
\begin{equation}\label{eq:growth-intro}
    \frac{1}{N}\sum_{i=1}^{N}\normf{\nabla f_i(W)}^2
    \leq D_1\normf{\nabla f(W)}^2+D_0,
    \qquad \forall W\in\R^{m\times n}.
\end{equation}
This growth condition includes bounded heterogeneity ($D_1=1$) and strong growth ($D_0=0$) as special cases and does not require bounded gradients~\cite{zhang2026adam}.
Under these assumptions, we establish both a positive and a negative result. For each fixed problem, Muon with $\beta$ sufficiently close to one and suitable stepsizes converges to stationary points when $D_0=0$, and to a neighborhood of them whose radius vanishes as $\beta\to1$ when $D_0>0$. In contrast, for each fixed $\beta$, Muon remains bounded away from stationarity on some problem with sufficiently many components.

We summarize our contributions as follows.
\begin{itemize}[leftmargin=1.5em]
    \item \textbf{Convergence for $\beta$ close to one.} If $\beta$ exceeds a problem-dependent threshold and the stepsize $\eta$ is constant, \cref{thm:recursion} gives
    \begin{equation}\label{eq:res_intro}
        \frac{1}{T}\sum_{t=0}^{T-1}\norm{\nabla f(W_t)}_* = \order{\frac{\Delta_f}{\eta T} + \frac{Lr\eta}{1-\beta} + (1-\beta)N\sqrt{rD_0} }
    \end{equation}
    for every sequence of permutations, up to terms that decay as $\orderi{1/T}$. Here, $\Delta_f=f(W_0)-f^*$, $L$ is the Lipschitz constant of the component gradients, and $r=\min\crki{m,n}$. 
    Choosing $1-\beta=\Theta\prti{(L\Delta_f/(N^2D_0T))^{1/3}}$ and $\eta=\Theta\prti{\sqrt{\Delta_f(1-\beta)/(LrT)}}$ yields an $\orderi{\varepsilon^{-3}}$ iteration complexity (\cref{cor:complexity}). To the best of our knowledge, this is the first such guarantee for Muon under shuffling, and its dependence on $\varepsilon$ matches that of SGD under shuffling (\Cref{tab:comparison}). 
    Under strong growth ($D_0=0$), suitable diminishing stepsizes and a single problem-dependent $\beta$ yield $\normi{\nabla f(W_t)}_*\to0$, without retuning $\beta$ when the target accuracy changes (\cref{cor:exact}).

    \item \textbf{Failure for every fixed $\beta$.} Let $N=4K+1$ for an integer $K\geq1$. If $\beta<(N+1)^{-1/K}$, then Muon fails on a smooth convex finite sum of $N$ components: for every sequence of permutations and every stepsize sequence satisfying $\eta_t\leq1/(2N)$ and remaining constant within each epoch, the iterates satisfy $\normi{\nabla f(W_t)}_*\geq1/N^2$ for all $t$. This problem has $1$-Lipschitz component gradients, satisfies \eqref{eq:growth-intro} with $D_1=1$ and $D_0\leq1$, and does not depend on $\beta$ (\cref{prop:smooth-n1}). As $N$ grows, the threshold $(N+1)^{-1/K}$ approaches one. Consequently, for any fixed $\beta<1$, Muon fails to reach stationarity on this family once $N$ is sufficiently large.

    \item \textbf{Implications for tuning $\beta$.} Our positive and negative results suggest two tuning rules: (i) if the training loss stops decreasing even after reducing the stepsize, increase $\beta$~(\cref{prop:smooth-n1}); (ii) when increasing $\beta$, decrease the stepsize accordingly~(\cref{cor:complexity}). \Cref{sec:exp-tuning} evaluates both suggestions.
\end{itemize}

\subsection{Related Work}

\textbf{Shuffling methods.}
Shuffling often accelerates SGD in practice~\cite{bottou2009curiously}. The
theoretical analysis of SGD under shuffling has progressed from establishing faster asymptotic rates on strongly convex finite-sum problems~\cite{gurbuzbalaban2021random} to deriving nonasymptotic iteration-complexity bounds for RR, SO, and IG in strongly convex, convex, and
nonconvex settings~\cite{tseng1998incremental,mishchenko2020random,nguyen2020unified,liu2026random}.
In particular, for smooth nonconvex finite-sum problems, SGD under shuffling attains an $\orderi{\varepsilon^{-3}}$ iteration complexity under growth conditions of the form \eqref{eq:growth-intro}~\cite{nguyen2020unified}, as well as under a more general growth condition in~\cite{mishchenko2020random}.
Shuffling has been studied in combination with momentum~\cite{tran2021smg,qiu2024random,liang2025last,xiao2026stochastic}, under the Kurdyka--\L{}ojasiewicz inequality~\cite{li2023convergence}, and in distributed settings~\cite{mishchenko2022proximal,huang2023drr,huang2023abcrr}. 

\textbf{Muon and spectral methods.}
Muon~\cite{jordan2024muon} can be viewed as steepest descent under the spectral norm~\cite{bernstein2024old}. This view connects Muon to Shampoo-type preconditioning~\cite{gupta2018shampoo}. Muon can also be interpreted as a non-Euclidean trust-region method~\cite{kovalev2025understanding} or through linear minimization oracles over norm balls~\cite{pmlr-v267-pethick25a,riabinin2025gluon}. \citet{ma2026schattor} place Muon within a broader family of Schatten-norm methods.
Beyond these interpretations, \citet{davis2025spectral} identify conditions under which spectral updates outperform Euclidean ones. Under sampling with replacement, exact-polar Muon with ordinary momentum attains an $\orderi{\varepsilon^{-4}}$ complexity~\cite{li2025note,shen2026on}, while momentum variance reduction improves this to $\orderi{\varepsilon^{-3}}$~\cite{qian2025muon}. The effects of using finitely many Newton--Schulz iterations and other polynomial approximations to the polar
map have been analyzed separately~\cite{kim2026convergence,do2026muon,wang2026convergence}. 

\textbf{Muon under RR.}
\citet{li2026implicit} prove that stochastic steepest descent under RR, including Muon, converges to an approximate max-margin solution for linear classification on separable data. \citet{kim2026geometry} show that, on Gaussian mixture data, row normalization achieves higher population accuracy than Muon under RR. Neither result establishes an $\orderi{\varepsilon^{-3}}$ iteration complexity for smooth nonconvex finite-sum problems.

\textbf{Counterexamples and problem-dependent tuning.}
Hyperparameters chosen before the problem is specified can cause adaptive methods to fail~\cite{reddi2018convergence}. Related failures are known for signSGD~\cite{karimireddy2019error}, for sign methods without momentum under RR~\cite{qin2023convergence}, and for Muon on convex Lipschitz problems~\cite{parshakova2026muon,lai2026nonsmooth}. \citet{lai2026nonsmooth} further show that, in the nonsmooth setting, asymptotic convergence can be restored by letting the momentum coefficient approach one while the stepsize vanishes. For Adam, \citet{zhang2026adam} show that such failures can be avoided when the hyperparameters are chosen with knowledge of the problem. Our results establish a similar conclusion for Muon under shuffling.

\section{Problem Setting and Algorithm}\label{sec:setting}

This section introduces the notation, states the assumptions, and specifies the algorithm.

\textbf{Notation.}
For $A,B\in\R^{m\times n}$, let $\iprodi{A}{B}=\tr(A^\T B)$ denote the Frobenius inner product. We use $\normif{A}$, $\normi{A}_*$, and $\normi{A}_2$ to denote the Frobenius, nuclear, and spectral norms, respectively. Let $r:=\min\crki{m,n}$. These norms satisfy $\normif{A}\leq\normi{A}_*\leq\sqrt{r}\normif{A}$ and $\normif{A}\leq\sqrt{r}\normi{A}_2$.

\textbf{Polar map.}
Let $M\in\R^{m\times n}$ be nonzero with rank $s$. Consider its compact singular value decomposition $M=U_s\Sigma_sV_s^\T$, where $U_s\in\R^{m\times s}$ and $V_s\in\R^{n\times s}$ have orthonormal columns, and $\Sigma_s\in\R^{s\times s}$ is diagonal with positive entries. The polar map is defined as
\begin{equation}\label{eq:polar}
    \Pol(M):=U_sV_s^\T,\quad \Pol(\mathbf{0}):=\mathbf{0}.
\end{equation}
This definition does not depend on the choice of the decomposition.

\textbf{Stationarity measure.}
As in prior analyses of Muon~\cite{shen2026on,li2025note}, we measure stationarity by $\normi{\nabla f(W)}_*$, which upper bounds $\normif{\nabla f(W)}$.

\subsection{Assumptions}\label{sec:assumptions}

We make two standard assumptions: componentwise smoothness~\cite{mishchenko2020random,nguyen2020unified,shen2026on} and a growth condition~\cite{bottou2018optimization,nguyen2020unified,zhang2026adam}.

\begin{assumption}\label[assumption]{ass:smooth}
Each $f_i:\R^{m\times n}\to\R$ is differentiable, and there exists $L>0$ such that, for all $W,Z$ and all $i$,
\begin{equation}\label{eq:smooth}
    \normf{\nabla f_i(W)-\nabla f_i(Z)}\leq L\normf{W-Z}.
\end{equation}
The average objective satisfies $f^*:=\inf_W f(W)>-\infty$.
\end{assumption}
\Cref{ass:smooth} implies that $f$ is $L$-smooth. Moreover, since $\norm{A}_*\leq\sqrt{r}\normif{A}$ and $\normif{B}\leq\sqrt{r}\norm{B}_2$, it also implies that each $f_i$ satisfies the spectral-norm smoothness condition of \citet{shen2026on} with constant $rL$:
$\normi{\nabla f_i(W)-\nabla f_i(Z)}_*\leq rL\normi{W-Z}_2$.

\begin{assumption}\label[assumption]{ass:growth}
There exist finite constants $D_1\geq0$ and $D_0\geq0$ such that
\begin{equation}
    \label{eq:growth}
    \frac{1}{N}\sum_{i=1}^N\normf{\nabla f_i(W)}^2
    \leq D_1\normf{\nabla f(W)}^2+D_0,\quad \forall W\in\R^{m\times n}.
\end{equation}
\end{assumption}
Define $\rho:=\max\crki{0,D_1-1}$. We have
\begin{equation}\label{eq:relative-variance}
    \frac{1}{N}\sum_{i=1}^N\normf{\nabla f_i(W)-\nabla f(W)}^2
    = \frac{1}{N}\sum_{i=1}^N\normf{\nabla f_i(W)}^2-\normf{\nabla f(W)}^2
    \leq \rho\normf{\nabla f(W)}^2+D_0.
\end{equation}
%At a stationary point of $f$, \eqref{eq:growth} reduces to $\frac{1}{N}\sum_{i=1}^N\normif{\nabla f_i(W)}^2\leq D_0$, so $D_0$ bounds the dispersion of the component gradients that remains at stationary points. 
When $D_1=1$, \eqref{eq:relative-variance} is a bounded heterogeneity condition. When
$D_0=0$, \eqref{eq:growth} is the strong growth
condition~\cite{vaswani2019fast}, under which every stationary point of $f$ is
a stationary point of every $f_i$.

\subsection{Exact-polar Muon under shuffling}\label{sec:algorithm}

\begin{algorithm}[htbp]
\caption{Exact-polar Muon under shuffling}\label{alg:rr-muon}
\begin{algorithmic}[1]
\Require $W_0\in\R^{m\times n}$, $M_{-1}=M_{\mathrm{init}}\in\R^{m\times n}$,
constant $\beta\in[0,1)$, positive stepsizes $\crki{\eta_t}_{t\geq0}$
\For{$k=0,1,2,\ldots$}
    \State Choose any permutation $\pi_k$ of $\crki{1,\ldots,N}$
    \For{$j=0,\ldots,N-1$}
        \State $t\gets kN+j$, \quad $i_t\gets\pi_k(j+1)$
        \State $G_t\gets\nabla f_{i_t}(W_t)$
        \State $M_t\gets\beta M_{t-1}+(1-\beta)G_t$
        \State $W_{t+1}\gets W_t-\eta_t\Pol(M_t)$
    \EndFor
\EndFor
\end{algorithmic}
\end{algorithm}

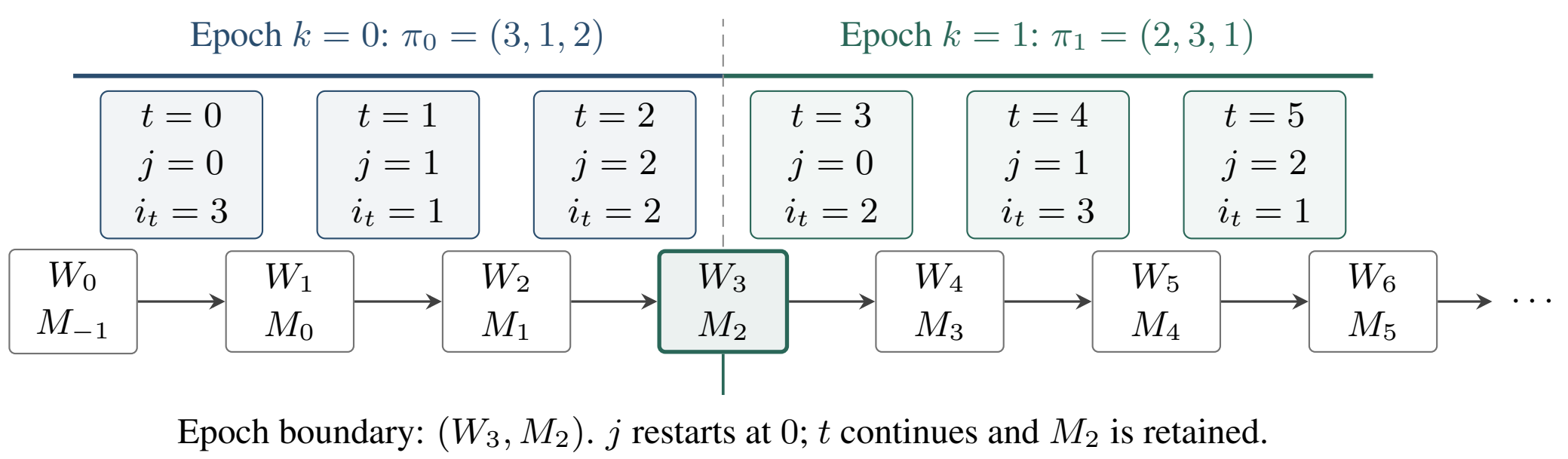
\begin{figure}[htbp]
    \centering
    % Vector illustration of Algorithm 1. The example has N=3.
\begingroup
\definecolor{rrepochzero}{RGB}{42,77,110}
\definecolor{rrepochone}{RGB}{42,103,88}
\begin{tikzpicture}[
    x=1.9cm,y=1cm,
    font=\small,
    >={Stealth[length=4pt,width=4pt]},
    update/.style={
        draw,rounded corners=2pt,align=center,
        minimum width=1.42cm,minimum height=1.18cm,
        inner sep=3pt,font=\footnotesize
    },
    state/.style={
        draw=black!55,fill=white,rounded corners=1.5pt,
        align=center,minimum width=1.12cm,minimum height=0.84cm,
        inner sep=3pt
    }
]
    % General indexing rule; the timeline below specializes to N=3.
    \node at (3,2.83)
        %{$t=kN+j,\qquad i_t=\pi_k(j+1)\qquad (N=3)$};
        {};
    % Epoch spans are aligned with the states at their endpoints.
    \draw[rrepochzero,line width=1.1pt] (0,1.98) -- (3,1.98);
    \draw[rrepochone,line width=1.1pt] (3,1.98) -- (6,1.98);
    \node[text=rrepochzero,anchor=south] at (1.5,2.04)
        {Epoch $k=0$: $\pi_0=(3,1,2)$};
    \node[text=rrepochone,anchor=south] at (4.5,2.04)
        {Epoch $k=1$: $\pi_1=(2,3,1)$};
    \draw[black!45,densely dashed] (3,0.48) -- (3,2.48);

    % Each upper card labels the update between two lower states.
    \foreach \x/\t/\j/\i in {0.5/0/0/3,1.5/1/1/1,2.5/2/2/2}{
        \node[update,draw=rrepochzero,fill=rrepochzero!6] at (\x,1.2)
            {$t=\t$\\[2pt]$j=\j$\\[2pt]$i_t=\i$};
    }
    \foreach \x/\t/\j/\i in {3.5/3/0/2,4.5/4/1/3,5.5/5/2/1}{
        \node[update,draw=rrepochone,fill=rrepochone!6] at (\x,1.2)
            {$t=\t$\\[2pt]$j=\j$\\[2pt]$i_t=\i$};
    }

    \node[align=right,anchor=east,font=\footnotesize] at (-0.46,0)
        %{State before\\the update\\$(W_t,M_{t-1})$};
        {};
    \foreach \t/\previous in {0/{-1},1/0,2/1,4/3,5/4,6/5}{
        \node[state] (s\t) at (\t,0)
            {$W_{\t}$\\[2pt]$M_{\previous}$};
    }
    \node[state,draw=rrepochone,line width=1pt,fill=rrepochone!9]
        (s3) at (3,0) {$W_3$\\[2pt]$M_2$};
    \foreach \t/\next in {0/1,1/2,2/3,3/4,4/5,5/6}{
        \draw[->,black!75,line width=0.6pt] (s\t.east) -- (s\next.west);
    }
    \draw[->,black!75,line width=0.6pt] (s6.east) -- (6.55,0);
    \node[anchor=west] at (6.57,0) {$\cdots$};

    \draw[rrepochone,line width=0.7pt] (s3.south) -- (3,-0.82);
    \node[anchor=north,align=center,font=\footnotesize] at (3,-0.88)
        {Epoch boundary: $(W_3,M_2)$.
        $j$ restarts at $0$; $t$ continues and $M_2$ is retained.};
\end{tikzpicture}
\endgroup
    \caption{Iteration indices and state transitions in \cref{alg:rr-muon}, illustrated for $N=3$ with two example
    permutations. The upper boxes label updates. The lower boxes show the parameter and momentum states before each update.}
    \label{fig:rr-muon-indexing}
\end{figure}

\Cref{alg:rr-muon} describes the method we analyze, which uses the exact polar map \eqref{eq:polar}. At the start of epoch $k$, a permutation $\pi_k$ of $\crki{1,\ldots,N}$ determines the order in which the components are visited.
At iteration $t=kN+j$, the algorithm evaluates the gradient of component $i_t=\pi_k(j+1)$ at $W_t$, incorporates it into the momentum buffer, and updates $W_t$ along $-\Pol(M_t)$ with stepsize $\eta_t$. The momentum buffer is initialized at an arbitrary $M_{\mathrm{init}}$ and carried across epochs. Since the sequence $\crki{\pi_k}_{k\geq0}$ is arbitrary, the analysis covers RR, SO, and IG. \Cref{fig:rr-muon-indexing} illustrates the indexing in \cref{alg:rr-muon} for $N=3$.

% \textbf{Relation to Practical Muon implementations}
Practical Muon implementations~\cite{jordan2024muon,pytorch2026muon} differ from \cref{alg:rr-muon} mainly in that they approximate the polar map using a few Newton--Schulz iterations. The other differences are configurable: in \texttt{torch.optim.Muon}, Nesterov momentum and weight decay can be disabled, and the shape-dependent factor multiplying the stepsize can be absorbed into $\eta_t$.
With these settings, each dataset item can be treated as one component, either a sample or a fixed mini-batch. Iterating over a \texttt{DataLoader} with \texttt{shuffle=True} and \texttt{batch\_size=1} then realizes \cref{alg:rr-muon} with RR, up to the Newton--Schulz approximation. \Cref{app:implementation} provides the full configuration.

\section{Main Results}
\label{sec:results}

In this section, we present the main results.
First, \cref{thm:recursion} gives an approximate one-step descent inequality and bounds the weighted average of $\normi{\nabla f(W_t)}_*$ in \eqref{eq:weighted-finite}, which yields the $\orderi{\varepsilon^{-3}}$ iteration complexity in \cref{cor:complexity}. We then discuss two choices of the momentum coefficient $\beta$: one that depends on the target accuracy (\cref{rem:tolerance}) and, under strong growth, one that does not (\cref{cor:exact}). Finally, \Cref{prop:smooth-n1} shows that $\beta$ cannot be chosen independently of the problem. Throughout, $\beta$ is held constant within each run, and all results hold for every sequence of permutations. Proofs are deferred to \cref{app:smooth-proofs,app:smooth-counterexample-proof}.

\subsection{Convergence result}\label{sec:convergence}

The bounds in the subsequent results depend on the problem through three constants:
\begin{equation}\label{eq:heterogeneity-constants}
    \Gamma_0:=N\sqrt{rD_0},\quad \Gamma_1:=N\sqrt{r\rho},\quad \Delta_f:=f(W_0)-f^*.
\end{equation}
Notably, these constants do not depend on $\beta$, $T$, or $\varepsilon$.

\begin{theorem}
\label{thm:recursion}
Let \Cref{ass:smooth,ass:growth} hold.
For any $\beta\in[0,1)$ and any positive, nonincreasing
stepsize sequence $\crki{\eta_t}$, the iterates $\crki{W_t}_{t\geq 0}$ generated by
\cref{alg:rr-muon} satisfy 
\begin{equation}\label{eq:descent}
\begin{aligned}
    f(W_{t+1}) &\leq f(W_t)
        -\eta_t\prt{1-\frac{2(1-\beta)\Gamma_1}{1-\beta^{t+1}}}
            \norm{\nabla f(W_t)}_*
        +\frac{2(1-\beta)\Gamma_0\eta_t}{1-\beta^{t+1}}\\
    &\quad +2Lr\eta_t
        \sum_{j=0}^{t-1}\beta^{t-j}\eta_j +\frac{2\sqrt{r}\eta_t\beta^{t+1}}{1-\beta^{t+1}}
        \normf{M_{\mathrm{init}}}+\frac{Lr}{2}\eta_t^2.
\end{aligned}
\end{equation}

%Define $\Delta_f:=f(W_0)-f^*.$
If $2(1-\beta)\Gamma_1 < 1$, then, for every $T\geq1$,
\begin{equation}\label{eq:weighted-finite}
\begin{aligned}
    \frac{\sum_{t=0}^{T-1}\eta_t\norm{\nabla f(W_t)}_*}
        {\sum_{t=0}^{T-1}\eta_t}
    &\leq\frac{\Delta_f+Lr\prt{\frac{1}{2}+\frac{2\beta}{1-\beta}}
            \sum_{t=0}^{T-1}\eta_t^2}
        {\prt{1-2(1-\beta)\Gamma_1}\sum_{t=0}^{T-1}\eta_t} + \frac{2\eta_0\beta
        \brk{\sqrt{r}\normf{M_{\mathrm{init}}}
            +(1-\beta)\prt{\Gamma_1\norm{\nabla f(W_0)}_*+\Gamma_0}}}
        {(1-\beta)^2\prt{1-2(1-\beta)\Gamma_1}
            \sum_{t=0}^{T-1}\eta_t}\\
    &\quad+\frac{2\Gamma_1 Lr\eta_0^2\beta^2}
        {(1-\beta)^2\prt{1-2(1-\beta)\Gamma_1}
            \sum_{t=0}^{T-1}\eta_t}
        +\frac{2(1-\beta)\Gamma_0}{1-2(1-\beta)\Gamma_1}.
\end{aligned}
\end{equation}

\end{theorem}

In \eqref{eq:weighted-finite}, both the condition $2(1-\beta)\Gamma_1<1$ and the final term
\begin{equation}\label{eq:residual}
    R:=\frac{2(1-\beta)\Gamma_0}{1-2(1-\beta)\Gamma_1}
\end{equation}
favor choosing $\beta$ close to one. In particular, $R$ does not decrease with the horizon $T$, but vanishes as $\beta\to1$. \Cref{cor:complexity} therefore lets $\beta$ approach one as the horizon $T$ grows and balances the remaining terms using a constant stepsize.

\begin{corollary}
\label[corollary]{cor:complexity}
Let \Cref{ass:smooth,ass:growth} hold, $M_{\mathrm{init}}=\mathbf{0}$, and $D_0>0$. For any integer $T\geq1$, assume $\Delta_f >0$ and choose
\begin{equation}\label{eq:horizon-parameters}
\begin{aligned}
    1-\beta&=\frac{1}{1+4\Gamma_1+  \prt{\frac{2\Gamma_0^2T}{Lr\Delta_f}}^{1/3}},\quad
    \eta_t=\eta:=\sqrt{\frac{\Delta_f(1-\beta)}{2LrT}},\qquad t\geq 0.
\end{aligned}
\end{equation}
Then
\begin{equation}\label{eq:horizon-rate-order}
\begin{aligned}
    \frac{1}{T}\sum_{t=0}^{T-1}\norm{\nabla f(W_t)}_*
    &\leq\mathcal{O}\left(\prt{\frac{Lr\Gamma_0\Delta_f}{T}}^{1/3}+\sqrt{\frac{(1+\Gamma_1)Lr\Delta_f}{T}}
        +\frac{\Gamma_0^{2/3}\prt{\Gamma_1\norm{\nabla f(W_0)}_*+\Gamma_0}}{(Lr\Delta_f)^{1/3}T^{2/3}}\right.\\
    &\quad\left.+\frac{(1+\Gamma_1)\prt{\Gamma_1\norm{\nabla f(W_0)}_*+\Gamma_0}}{T}
        +\frac{\Gamma_1(1+\Gamma_1)^{3/2}\sqrt{Lr\Delta_f}}{T^{3/2}}\right),
\end{aligned}
\end{equation}
where $\orderi{\cdot}$ hides numerical constants. Moreover, the number of iterations required to achieve $\frac{1}{T}\sum_{t=0}^{T-1}\normi{\nabla f(W_t)}_*\leq\varepsilon$ is
\begin{equation}\label{eq:complexity-horizon}
\begin{aligned}
   &\mathcal{O}\left(\frac{Lr\Gamma_0\Delta_f}{\varepsilon^3}+\frac{(1+\Gamma_1)Lr\Delta_f}{\varepsilon^2}
        +\frac{\Gamma_0\prt{\Gamma_1\norm{\nabla f(W_0)}_*+\Gamma_0}^{3/2}}{\sqrt{Lr\Delta_f}\,\varepsilon^{3/2}}\right.\\
   &\quad\left.+\frac{(1+\Gamma_1)\prt{\Gamma_1\norm{\nabla f(W_0)}_*+\Gamma_0}}{\varepsilon}
        +\frac{\Gamma_1^{2/3}(1+\Gamma_1)(Lr\Delta_f)^{1/3}}{\varepsilon^{2/3}}\right).
\end{aligned}
\end{equation}
\end{corollary}

The $\orderi{\varepsilon^{-3}}$ iteration complexity matches that of SGD under shuffling for smooth nonconvex finite-sum problems~\cite{mishchenko2020random,nguyen2020unified} in its dependence on $\varepsilon$.

\begin{remark}
\label[remark]{rem:tolerance}
In \cref{cor:complexity}, $1-\beta=\Theta\prti{(L\Delta_f/(N^2D_0T))^{1/3}}$ and $\eta=\Theta\prti{\sqrt{\Delta_f(1-\beta)/(LrT)}}$ depend on the horizon $T$, and hence
on the target accuracy: $1-\beta=\Theta(\varepsilon)$ and $\eta=\Theta(\varepsilon^2)$. Thus, achieving a smaller target accuracy requires retuning both parameters. The choice of $\eta_t$ in \eqref{eq:horizon-parameters} motivates tuning suggestion (ii): when increasing $\beta$, decrease the stepsize accordingly.
\end{remark}

In contrast, under strong growth ($D_0=0$), a single constant $\beta$ suffices for all target accuracies (\cref{cor:exact}).

\begin{corollary}
\label[corollary]{cor:exact}
In the setting of \cref{thm:recursion}, suppose $D_0=0$, $M_{\mathrm{init}} = \mathbf{0}$ and
\begin{equation}\label{eq:steps-exact}
    \sum_{t=0}^\infty\eta_t=\infty,
    \qquad \sum_{t=0}^\infty\eta_t^2<\infty.
\end{equation}
If $D_1>1$, choose a momentum coefficient satisfying
\begin{equation}\label{eq:beta-strict}
    1>\beta>1-\frac{1}{2\Gamma_1}.
\end{equation}
If $D_1\leq1$, any constant $\beta\in[0,1)$ is admissible.
Then
\begin{equation}\label{eq:exact-conclusion}
    \norm{\nabla f(W_t)}_*\to0,
    \qquad f(W_t)\to f_\infty\in[f^*,\infty).
\end{equation}
Every accumulation point of $\crki{W_t}_{t\geq 0}$ is stationary. These conclusions hold for every sequence of permutations.
\end{corollary}

\subsection{Failure for a prescribed coefficient}\label{sec:counterexample}

The choices of $\beta$ in \cref{cor:complexity,cor:exact} are made after the problem is specified, since they depend on $N$. We now show that this order cannot be reversed: for every prescribed $\beta\in[0,1)$, Muon remains bounded away from stationarity on a family of smooth convex finite-sum problems over $2\times2$ matrices.

\begin{theorem}
\label[theorem]{prop:smooth-n1}
Fix $\beta\in[0,1)$, and choose an integer $K\geq1$ such that $(4K+2)\beta^K<1$. Let $N=4K+1$.
Define the Huber loss $h(u)=u^2/2$ for $|u|\leq1$ and $h(u)=|u|-1/2$ otherwise, and let
\begin{equation}\label{eq:huber-components}
    g:=\frac{1}{N}\sum_{i=1}^{N}g_i,\quad
    g_1(x):=N h(x+2),\quad
    g_i(x):=h(x-2),\qquad2\leq i\leq N.
\end{equation}
Define
\begin{equation}\label{eq:huber-matrix}
    f_i(W):=g_i([W]_{11})
        +\half\sum_{(a,b)\neq(1,1)}[W]_{ab}^2,
    \quad f:=\frac{1}{N}\sum_{i=1}^{N}f_i,\qquad W\in\R^{2\times2},
\end{equation}
where $[W]_{ab}$ denotes the $(a,b)$ entry of $W$. Then the following hold:
\begin{enumerate}
\item[(i)] Each $f_i$ is convex with an $N$-Lipschitz gradient, and $f_1,\ldots,f_N$ satisfy \Cref{ass:growth} with $D_1=1$ and
$D_0=(N-1)(N+1)^2/N^2$.
\item[(ii)] Run \cref{alg:rr-muon} from $W_0=M_{-1}=0$, using any stepsizes
that are constant within each epoch and satisfy
\begin{equation}\label{eq:huber-stepsizes}
    \eta_{kN+j}=\eta_{kN}\in\Bigl(0,\frac{1}{2N}\Bigr],\qquad 0\leq j<N,\quad k\geq0.
\end{equation}
Then, for every sequence of permutations and every $t\geq0$,
\begin{equation}\label{eq:huber-failure}
    \norm{\nabla f(W_t)}_*\geq\frac{1}{N},\qquad
    f(W_t)-f^*\geq\frac{N+1}{2N^2}>0.
\end{equation}
\item[(iii)] The rescaled components $\tilde{f}_i:=f_i/N$ have $1$-Lipschitz gradients and satisfy \Cref{ass:growth} with $D_1=1$ and $D_0=(N-1)(N+1)^2/N^4\leq1$. Under the conditions in (ii), \cref{alg:rr-muon} generates the same iterates on $\tilde{f}_1,\ldots,\tilde{f}_N$. Thus, for $\tilde{f}:=f/N$, we have $\normi{\nabla\tilde{f}(W_t)}_*\geq1/N^2$ for all $t\geq0$.
\end{enumerate}
\end{theorem}

\begin{remark}\label[remark]{rem:fixed-constants}
The condition $(4K+2)\beta^K<1$ in \cref{prop:smooth-n1} is equivalent to $\beta<(N+1)^{-4/(N-1)}$, whose right-hand side approaches one as $N$ grows. For the rescaled family of \cref{prop:smooth-n1}(iii), \Cref{ass:smooth,ass:growth} hold with $L=1$, $D_1=1$, and $D_0=1$ for every $N$. Hence, no fixed $\beta<1$ works for the entire family, and $\beta$ must depend on $N$.
\end{remark}

\begin{remark}\label[remark]{rem:stepsizes}
Condition \eqref{eq:huber-stepsizes} includes every constant stepsize
$\eta\in(0,1/(2N)]$ and diminishing stepsizes such as
\begin{equation}\label{eq:huber-steps}
    \eta_{kN+j}=\frac{1}{2N(k+1)^p},\qquad 0\leq j<N,\quad \frac{1}{2}<p\leq1,
\end{equation}
which satisfy $\sum_t\eta_t=\infty$ and $\sum_t\eta_t^2<\infty$. Thus, neither small nor diminishing stepsizes eliminate the failure. This motivates tuning suggestion (i): if the training loss stops decreasing even after the stepsize is reduced, increase $\beta$.
\end{remark}

\section{Experimental Results}\label{sec:experiments}

We report three groups of experiments. \Cref{sec:exp-orders} compares sampling schemes.
\Cref{sec:exp-batch} examines how the required momentum coefficient depends on the number of components, and
\Cref{sec:exp-tuning} evaluates the two tuning suggestions from \cref{sec:intro,sec:convergence}.
Unless stated otherwise, we use RR. Experimental details are provided in \cref{app:experiments}.

\subsection{Sampling schemes}\label{sec:exp-orders}

We compare RR, SO, IG, and sampling with replacement on three models.
For each model, we partition the training data once into $N$ fixed batches of
size $B$, with $f_i$ being the average loss over the $i$-th batch.
All four schemes use the same components and constant stepsize, differing only in how the components are sampled.

\textbf{MLP on MNIST.} We train a two-layer tanh multilayer perceptron (MLP) on the $S=60000$ MNIST training images~\cite{mnist}, with an $\ell_2$ term in every component, using the exact polar map. For $B=100$ ($N=600$), we consider two data partitions: a random partition and a sorted partition, obtained by sorting the images by label. \Cref{tab:exp-mnist-orders} reports the final training loss and the gradient norm. Across all four settings, corresponding to the two partitions and two values of $\beta$, RR and SO achieve lower training loss and gradient norm than sampling with replacement. Increasing $\beta$ from $0.9$ to $0.99$ reduces both metrics for every sampling scheme. Notably, IG performs competitively with the random partition but degrades substantially with the sorted partition. These results suggest that random permutations can mitigate the adverse effects of an unfavorable data order, while a larger momentum coefficient can further reduce both metrics.

\begin{table}[!htbp]
    \centering
    \small
    \caption{Final training loss $f(W_T)$ and gradient norm $\normi{\nabla f(W_T)}_*$ for an MLP on MNIST. Values are medians over three random seeds; IG is deterministic and was run
    once. The lowest value in each column is shown in bold.}
    \label{tab:exp-mnist-orders}
    \begin{tabular}{@{}lcccccccc@{}}
        \toprule
        & \multicolumn{4}{c}{random partition}
        & \multicolumn{4}{c}{sorted partition} \\
        \cmidrule(lr){2-5}\cmidrule(l){6-9}
        & \multicolumn{2}{c}{$\beta=0.9$} & \multicolumn{2}{c}{$\beta=0.99$}
        & \multicolumn{2}{c}{$\beta=0.9$} & \multicolumn{2}{c}{$\beta=0.99$} \\
        \cmidrule(lr){2-3}\cmidrule(lr){4-5}\cmidrule(lr){6-7}\cmidrule(l){8-9}
        & $f$ & $\normi{\nabla f}_*$ & $f$ & $\normi{\nabla f}_*$
        & $f$ & $\normi{\nabla f}_*$ & $f$ & $\normi{\nabla f}_*$ \\
        \midrule
        RR & $\mathbf{0.1040}$ & $0.087$ & $0.0601$ & $0.061$ & $\mathbf{0.1447}$ & $\mathbf{0.149}$ & $0.0671$ & $0.087$ \\
        SO & $0.1042$ & $0.092$ & $0.0601$ & $0.070$ & $0.1468$ & $0.156$ & $\mathbf{0.0668}$ & $\mathbf{0.079}$ \\
        IG & $0.1050$ & $\mathbf{0.083}$ & $\mathbf{0.0600}$ & $\mathbf{0.059}$ & $2.2654$ & $1.154$ & $0.3523$ & $0.856$ \\
        with replacement & $0.1049$ & $0.099$ & $0.0638$ & $0.079$ & $0.1475$ & $0.176$ & $0.0711$ & $0.092$ \\
        \bottomrule
    \end{tabular}
\end{table}

\begin{figure}[htbp]
    \centering
    \begin{minipage}[b]{0.49\linewidth}
        \centering
        \includegraphics[width=\linewidth]{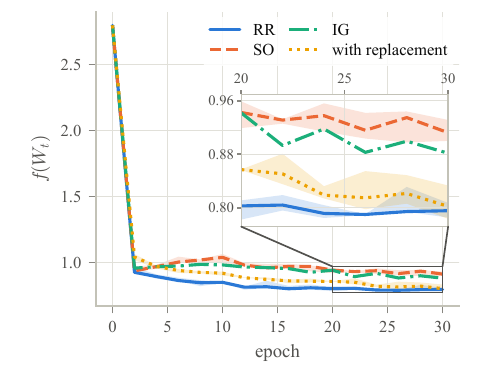}
    \end{minipage}\hfill
    \begin{minipage}[b]{0.49\linewidth}
        \centering
        \includegraphics[width=\linewidth]{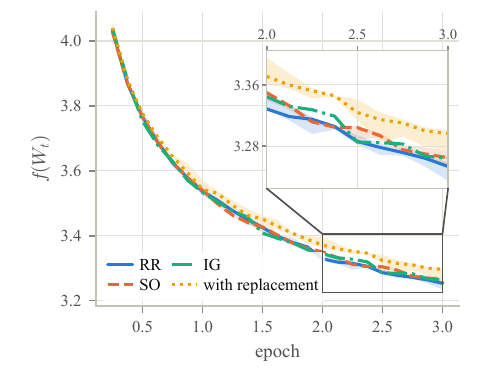}
    \end{minipage}
    \caption{Training loss $f(W_t)$ under four sampling schemes with $\beta=0.95$. Curves show the median over three seeds, and shaded bands show the range across seeds. Left: ResNet-18 on CIFAR-10 with $B=125$ ($N=400$). Right: GPT-2 small on FineWeb-Edu with $B=32$ ($N=9155$). }
    \label{fig:exp-cifar-gpt2}
\end{figure}

\textbf{ResNet-18 on CIFAR-10.} We train a ResNet-18~\cite{he2016deep} on CIFAR-10~\cite{krizhevsky2009learning} ($S=50000$ images) for $30$ epochs, using $B=125$ ($N=400$). All runs use $\beta=0.95$ and the Newton--Schulz approximation to the polar map. RR attains the lowest training loss throughout training (\Cref{fig:exp-cifar-gpt2}, left). After the first few epochs, SO and IG have higher training losses than sampling with replacement.

\textbf{GPT-2 on FineWeb-Edu.} We train GPT-2 small~\cite{radford2019language} for three epochs on $3\times 10^8$ tokens from FineWeb-Edu~\cite{penedo2024fineweb}, divided into sequences of $1024$ tokens, with $B=32$ sequences per batch ($N=9155$). In the later epochs, RR, SO, and IG have similar training losses, all below that of sampling with replacement. RR achieves the lowest final training loss (\Cref{fig:exp-cifar-gpt2}, right).

\subsection{Required momentum and the number of components}\label{sec:exp-batch}

\Cref{cor:complexity} selects $\beta$ closer to one as $N$ grows, while \cref{prop:smooth-n1} shows that every fixed $\beta$ fails on some problem once $N$ is sufficiently large. We examine this dependence using the counterexample from \eqref{eq:huber-matrix} and an MLP trained on MNIST.

\begin{figure}[!htbp]
    \centering
    \begin{minipage}[b]{0.49\linewidth}
        \centering
        \includegraphics[width=\linewidth]{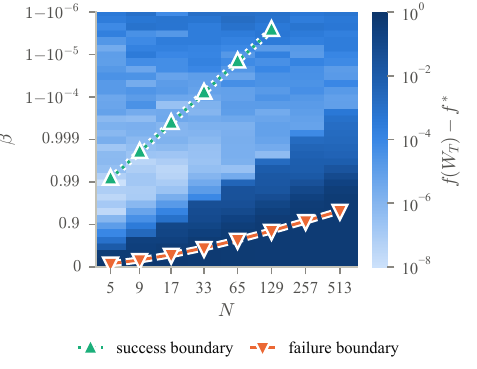}
    \end{minipage}\hfill
    \begin{minipage}[b]{0.49\linewidth}
        \centering
        \includegraphics[width=\linewidth]{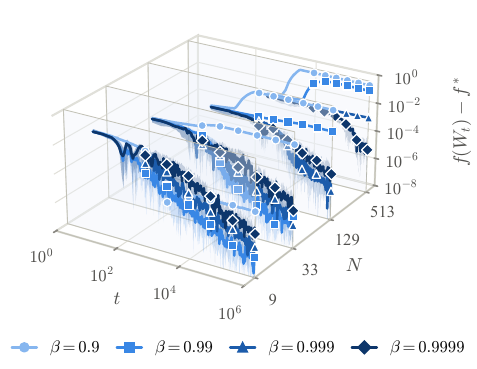}
    \end{minipage}
    \caption{The counterexample in \eqref{eq:huber-matrix} under RR.
    Left: final objective gap $f(W_T)-f^*$, shown as the median over
    five seeds. Right: objective gap $f(W_t)-f^*$ over the course of training for $N=9,33,129,513$ and four momentum coefficients, shown as the median over ten seeds.}
    \label{fig:exp-counterexample}
\end{figure}

\textbf{Smooth convex counterexample.} We run Muon under RR on the counterexample from \eqref{eq:huber-matrix} for $N\in\crki{5,9,17,\ldots,513}$, varying $1-\beta$ from $10^{-6}$ to $1$. All runs use the stepsizes in \eqref{eq:huber-steps} with $p=0.6$ and $T=10^6$ iterations.
\Cref{fig:exp-counterexample} (left) shows the final gap $f(W_T)-f^*$, together with two theoretical boundaries. Below the failure boundary ($(N+1)\beta^{(N-1)/4}=1$), \cref{prop:smooth-n1} guarantees that $\norm{\nabla f(W_t)}_*\geq1/N$ for all $t$. Above the success boundary ($2N(N+1)\sqrt{N-1}(1-\beta)=1$), the residual $R$ in \eqref{eq:residual} is below $1/N$. The observed transition from failure to convergence lies between these two boundaries and shifts toward $\beta=1$ as $N$ increases.
The trajectories in \Cref{fig:exp-counterexample} (right) show the same
transition.

\textbf{MLP on MNIST.} We train the MLP from \cref{sec:exp-orders} using Muon with RR on a random partition of the images into $N=S/B$ fixed batches of size $B\in\crki{25,50,100,200,400,800}$, giving $N\in\crki{2400,1200,600,300,150,75}$. For every batch size, the objective $f$ is the average loss over all $S$ images, so only the number of components changes. All runs use a stepsize of $\eta=3\times 10^{-4}$ and $T=2\times 10^4$ iterations. 
%so a smaller $N$ corresponds to more epochs.

\begin{figure}[!htbp]
    \centering
    \begin{minipage}[b]{0.49\linewidth}
        \centering
        \includegraphics[width=\linewidth]{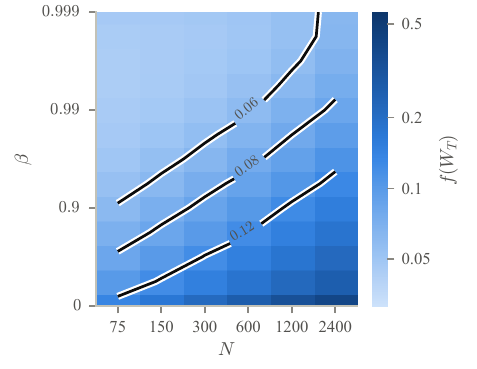}
    \end{minipage}\hfill
    \begin{minipage}[b]{0.49\linewidth}
        \centering
        \includegraphics[width=\linewidth]{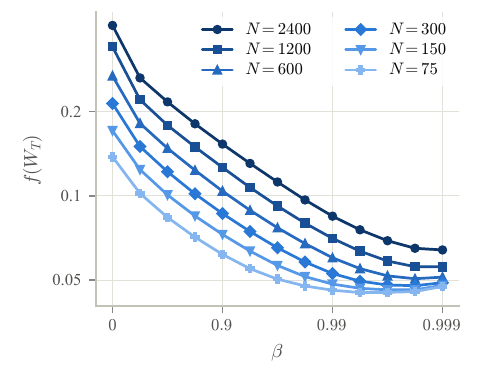}
    \end{minipage}
    \caption{Final training loss $f(W_T)$ for the MLP on MNIST under RR, with $\eta=3\times10^{-4}$ and $T=2\times10^4$. Values are medians over three seeds. Left: heatmap over $N=S/B$ and $\beta$, with empirical level curves. Right: the same values plotted against $\beta$ for each $N$.}
    \label{fig:exp-batch}
\end{figure}

\Cref{fig:exp-batch} (left) shows the final training loss over $(N,\beta)$, with level curves at $f(W_T)=0.06$, $0.08$, and $0.12$. These curves rise with $N$: achieving the same training loss requires a larger $\beta$ as $N$ increases. Thus, as in the counterexample \Cref{fig:exp-counterexample}, the MLP exhibits the same dependence: the $\beta$ required to achieve a given accuracy approaches one as the number of components grows.
\Cref{fig:exp-batch} (right) further illustrates this relationship: for fixed $N$, larger $\beta$ yields lower final training loss, whereas for fixed $\beta$, the final training loss increases with $N$.

\subsection{Tuning suggestions}\label{sec:exp-tuning}

We test suggestion (i) on the MLP, and suggestion (ii) on both the MLP and GPT-2 small.

\begin{figure}[!htbp]
    \centering
    \begin{minipage}[b]{0.49\linewidth}
        \centering
        \includegraphics[width=\linewidth]{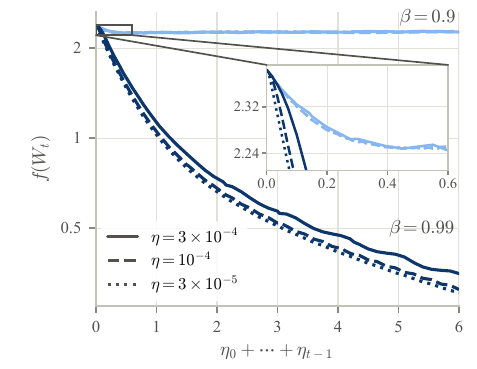}
    \end{minipage}\hfill
    \begin{minipage}[b]{0.49\linewidth}
        \centering
        \includegraphics[width=\linewidth]{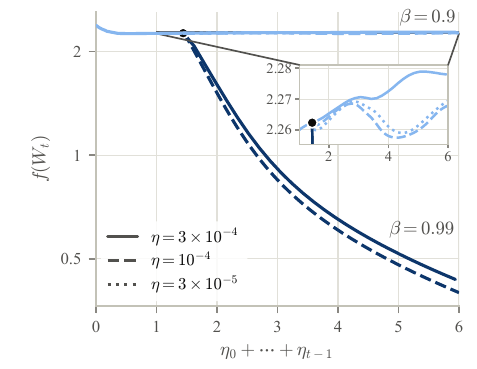}
    \end{minipage}
    \caption{Training loss $f(W_t)$ for the MLP on MNIST under IG with the sorted partition ($B=100$), plotted against the cumulative stepsize. Left: separate runs for each $(\beta,\eta)$. Right: the run with $\beta=0.9$ and $\eta=3\times10^{-4}$, continued from the plateau (black dot) using each $(\beta,\eta)$ pair.}
    \label{fig:exp-tuning-budget}
\end{figure}

\textbf{Suggestion (i): If the training loss stops decreasing even after reducing the stepsize, increase $\beta$.} We return to IG on the sorted partition with $B = 100$ ($N=600$), where $\beta=0.9$ yields a high training loss (\cref{tab:exp-mnist-orders}). We reduce the stepsize from $3\times10^{-4}$ to $10^{-4}$ and $3\times10^{-5}$, increasing $T$ proportionally, so that the total stepsize $\eta T=6$ is the same. With $\beta=0.9$, the training loss quickly plateaus and remains there for every stepsize. With $\beta = 0.99$, the loss continues to decrease (\Cref{fig:exp-tuning-budget}, left).
\Cref{fig:exp-tuning-budget} (right) further illustrates this behavior. Starting from the plateau at $\beta=0.9$ and $\eta = 3\times 10^{-4}$, we change $\eta$ and $\beta$. Increasing $\beta$ to $0.99$ substantially lowers the training loss, whereas reducing the stepsize alone does not lower the plateau.

\begin{figure}[!htbp]
    \centering
    \begin{minipage}[b]{0.49\linewidth}
        \centering
        \includegraphics[width=\linewidth]{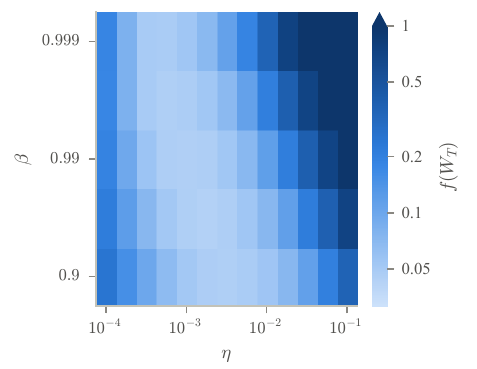}
    \end{minipage}\hfill
    \begin{minipage}[b]{0.49\linewidth}
        \centering
        \includegraphics[width=\linewidth]{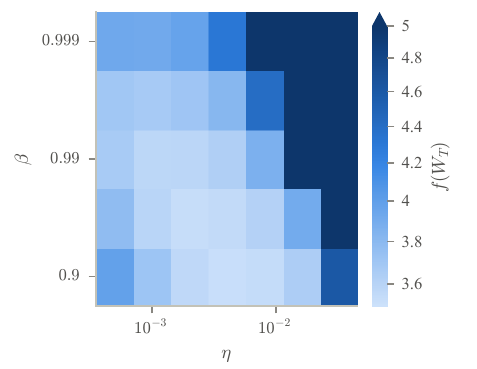}
    \end{minipage}
    \caption{Final training loss $f(W_T)$ under RR across values of $(\beta,\eta)$.
    Left: MNIST with $B=100$ and a random partition,
    $T=2\times 10^4$; values are medians over three seeds. Right: GPT-2 small trained for one epoch ($T=9155$); results are from one seed.}
    \label{fig:exp-tuning}
\end{figure}

\textbf{Suggestion (ii): When increasing $\beta$, decrease the stepsize as well.} 
We screen the final training loss $f(W_T)$ over a grid of $(\beta,\eta)$ pairs for MLP on MNIST ($B=100$, random partition) and for GPT-2 small. \Cref{fig:exp-tuning} shows that the stepsize yielding the lowest final training loss decreases as $\beta$ approaches one. This observation aligns with the coupling $\eta\propto\sqrt{1-\beta}$ in \eqref{eq:horizon-parameters}.

\FloatBarrier

\section{Conclusion}\label{sec:discussion}

We establish an $\orderi{\varepsilon^{-3}}$ iteration complexity for exact-polar Muon under shuffling on smooth nonconvex finite-sum problems. To the best of our knowledge, this is the first complexity guarantee for Muon in this setting, and its dependence on $\varepsilon$ matches that of SGD under shuffling. The guarantee requires a problem-dependent $\beta$ that approaches one as the number of components grows. In contrast, for any fixed $\beta\in[0,1)$, Muon remains bounded away from stationarity on a smooth convex finite-sum problem with sufficiently many components, regardless of the permutation sequence and even with small stepsizes. These results suggest two tuning rules: (i) if the training loss stops decreasing even after reducing the stepsize, increase $\beta$; and (ii) when increasing $\beta$, reduce the stepsize as well. Our experiments support both recommendations.

\FloatBarrier

\section*{Acknowledgments}

Computational resources were provided by the National Academic Infrastructure for Supercomputing in Sweden (NAISS), funded by the Swedish Research Council.

\appendix
\section{Proofs of the Convergence Results}\label[appendix]{app:smooth-proofs}

For any fixed sequence of permutations, all arguments in this appendix are deterministic.

\subsection{Proof of Theorem \ref{thm:recursion}}
\begin{lemma}
    \label{lem:polar_map}
    For every $M\in\R^{m\times n}$, the polar map in \eqref{eq:polar} satisfies
\begin{equation}\label{eq:polar-properties}
    \normf{\Pol(M)}\leq\sqrt{r},\qquad \norm{\Pol(M)}_*\leq r,\quad
    \iprod{M}{\Pol(M)}=\norm{M}_*,\qquad
    \Pol(aM)=\Pol(M)\quad(a>0).
\end{equation}
\end{lemma}

\begin{proof}
   Let $M\neq\mathbf{0}$ have the compact SVD $M = U_s\Sigma_sV_s^\T$ with $s\leq r$. Then $\Pol(M)=U_sV_s^\T$ has $s$ singular values, all equal to one, so
\begin{equation*}
    \normf{\Pol(M)} = \sqrt{\tr(U_sV_s^\T V_sU_s^\T)} = \sqrt{\tr(I_s)} = \sqrt{s} \leq \sqrt{r},\qquad
    \norm{\Pol(M)}_* = \tr(\sqrt{V_sU_s^\T U_sV_s^\T}) = \tr(I_s) = s \leq r.
\end{equation*}
Moreover, $\iprod{M}{U_sV_s^\T}=\tr(V_s\Sigma_sV_s^\T)=\tr(\Sigma_s)=\norm{M}_*$.
For $a>0$, $aM = U_s(a\Sigma_s)V_s^\T$ is the compact SVD of $aM$. Hence, $\Pol(aM) = \Pol(M)$. All four
relations hold trivially for $M=\mathbf{0}$.
\end{proof}

The proof of \cref{thm:recursion} begins with a one-step descent inequality, stated below.

\begin{lemma}\label[lemma]{lem:descent}
Let \Cref{ass:smooth} hold and $\beta\in[0,1)$. Then, the iterates $\crki{W_t}$ generated by \cref{alg:rr-muon} satisfy
\begin{equation}\label{eq:descent-decomp}
\begin{aligned}
    f(W_{t+1}) &\leq f(W_t)-\eta_t\norm{\nabla f(W_t)}_*+\frac{2\sqrt{r} \eta_t}{1-\beta^{t+1}}
        \normf{(1-\beta)\sum_{s=0}^t\beta^{t-s}
            \brk{\nabla f_{i_s}(W_t)-\nabla f(W_t)}}\\
    &\quad +\frac{2\sqrt{r} \eta_t(1-\beta)}{1-\beta^{t+1}}
        \sum_{s=0}^t\beta^{t-s}
            \normf{\nabla f_{i_s}(W_s)-\nabla f_{i_s}(W_t)}+\frac{2\sqrt{r} \eta_t\beta^{t+1}}{1-\beta^{t+1}}
        \normf{M_{\mathrm{init}}}+\frac{Lr}{2}\eta_t^2.
\end{aligned}
\end{equation}
\end{lemma}

\begin{proof}

    By \eqref{eq:polar-properties} and the update rule, we have
    \begin{equation}\label{eq:step-length}
        \normf{W_{t+1}-W_t} = \normf{-\eta_t \Pol(M_t)}\leq\sqrt{r}\eta_t.
    \end{equation}

    Applying the descent lemma to $W_{t + 1} = W_t - \eta_t \Pol(M_t)$ gives
    \begin{equation}
        \label{eq:descent_s1}
        \begin{aligned}
            f(W_{t+1}) &\leq f(W_t) - \eta_t\iprod{\nabla f(W_t)}{\Pol(M_t)} +\frac{Lr}{2}\eta_t^2.
        \end{aligned}
    \end{equation}

We now estimate the inner product in \eqref{eq:descent_s1}. Since $1-\beta^{t+1}>0$, it follows from \eqref{eq:polar-properties} that
\begin{equation}\label{eq:alignment}
\begin{aligned}
    \iprod{\nabla f(W_t)}{\Pol(M_t)}
    &= \iprod{\frac{M_t}{1-\beta^{t+1}}}{\Pol(M_t)}
        +\iprod{\nabla f(W_t)-\frac{M_t}{1-\beta^{t+1}}}{\Pol(M_t)}\\
    &= \norm{\frac{M_t}{1-\beta^{t+1}}}_*
        +\iprod{\nabla f(W_t)-\frac{M_t}{1-\beta^{t+1}}}{\Pol(M_t)}\\
    &\geq \norm{\frac{M_t}{1-\beta^{t+1}}}_*
        -\normf{\frac{M_t}{1-\beta^{t+1}}-\nabla f(W_t)}
            \normf{\Pol(M_t)}\\
    &\geq \norm{\frac{M_t}{1-\beta^{t+1}}}_*
        -\sqrt{r}\normf{\frac{M_t}{1-\beta^{t+1}}-\nabla f(W_t)}\\
    &\geq \norm{\nabla f(W_t)}_*
        -\norm{\frac{M_t}{1-\beta^{t+1}}-\nabla f(W_t)}_*
        -\sqrt{r}\normf{\frac{M_t}{1-\beta^{t+1}}-\nabla f(W_t)}\\
    &\geq \norm{\nabla f(W_t)}_*
        -2\sqrt{r}\normf{\frac{M_t}{1-\beta^{t+1}}-\nabla f(W_t)}.
\end{aligned}
\end{equation}
The first inequality is the Cauchy--Schwarz inequality and the second uses $\normif{\Pol(M_t)}\leq\sqrt{r}$. The third follows from the triangle inequality for the nuclear norm, and the last uses $\normi{A}_*\leq\sqrt{r}\normif{A}$. Expanding $M_t = \beta M_{t-1} + (1-\beta)\nabla f_{i_t}(W_t)$ recursively and noting $M_{-1}=M_{\mathrm{init}}$ gives
\begin{equation}\label{eq:tracking-decomp}
\begin{aligned}
    &\frac{M_t}{1-\beta^{t+1}}-\nabla f(W_t)= \frac{\beta^{t+1}M_{-1}}{1-\beta^{t+1}}
        +\frac{1-\beta}{1-\beta^{t+1}}\sum_{s=0}^t
            \beta^{t-s}\nabla f_{i_s}(W_s)-\nabla f(W_t)\\
    &= \frac{\beta^{t+1}M_{\mathrm{init}}}{1-\beta^{t+1}}+\frac{1-\beta}{1-\beta^{t+1}}\sum_{s=0}^t\beta^{t-s}\brk{\nabla f_{i_s}(W_s)-\nabla f_{i_s}(W_t)} 
    +\frac{1-\beta}{1-\beta^{t+1}}\sum_{s=0}^t\beta^{t-s}\brk{\nabla f_{i_s}(W_t)-\nabla f(W_t)}.
\end{aligned}
\end{equation}
Applying the triangle inequality to \eqref{eq:tracking-decomp} and
substituting the resulting inequality, together with \eqref{eq:alignment}, into \eqref{eq:descent_s1} yields \eqref{eq:descent-decomp}.

\end{proof}

\Cref{lem:rr-descent} bounds the third term on the right-hand side of \eqref{eq:descent-decomp} using \cref{lem:rr}. \cref{lem:rr} applies because each epoch visits every component exactly once.

\begin{lemma}\label[lemma]{lem:rr-descent}
Let \Cref{ass:smooth,ass:growth} hold and $\beta\in[0,1)$. Then, the iterates $\crki{W_t}_{t\geq 0}$  generated by \Cref{alg:rr-muon} satisfy
\begin{equation}\label{eq:rr-descent-bound}
\begin{aligned}
    &\frac{2\sqrt{r}\eta_t}{1-\beta^{t+1}}
        \normf{(1-\beta)\sum_{s=0}^t\beta^{t-s}
            \brk{\nabla f_{i_s}(W_t)-\nabla f(W_t)}}\leq\frac{2(1-\beta)\eta_t}{1-\beta^{t+1}}
        \prt{\Gamma_1\norm{\nabla f(W_t)}_*+\Gamma_0},
\end{aligned}
\end{equation}
where $\Gamma_0 = N\sqrt{rD_0}$ and $\Gamma_1 = N\sqrt{r\rho}$.
\end{lemma}
\begin{proof}
Fix $W\in\R^{m\times n}$.
For every $j\geq0$, we can write $j + 1 = qN + \ell$ with $q\geq0$ and $0\leq \ell<N$. The first $qN$ indices $s = 0,1,\ldots, qN-1$ cover the complete epochs. Hence, $\sum_{s=0}^{qN-1}\brki{\nabla f_{i_s}(W) - \nabla f(W)} = \mathbf{0}$. The remaining $i_{qN}, i_{qN+1}, \ldots, i_{qN+\ell-1}$ form a subset of $\crki{1,\ldots,N}$ without repetition. Thus, \cref{lem:rr} yields
\begin{equation}
    \label{eq:sumj}
    \normf{\sum_{s=0}^j
        \brk{\nabla f_{i_s}(W)-\nabla f(W)}}
    \leq \frac{1}{2}\sqrt{N\sum_{i=1}^N
        \normf{\nabla f_i(W)-\nabla f(W)}^2}.
\end{equation}
However, the weight $\beta^{t-s}$ prevents a complete epoch from canceling. We therefore express the weighted sum in terms of the unweighted prefix sums.
For $0<\beta<1$ and $0\leq s\leq t$, with empty sums interpreted
as zero, we have
\begin{equation}
    \label{eq:beta_ts}
    1-(1-\beta)\sum_{j=s}^{t-1}\beta^{t-1-j}
    =1-(1-\beta^{t-s})=\beta^{t-s}.
\end{equation}

In light of \eqref{eq:beta_ts}, we obtain
\begin{equation}\label{eq:weight_decomp}
\begin{aligned}
    &\sum_{s=0}^t\beta^{t-s}
        \brk{\nabla f_{i_s}(W)-\nabla f(W)}= \sum_{s=0}^t
        \brk{1-(1-\beta)\sum_{j=s}^{t-1}\beta^{t-1-j}}
        \brk{\nabla f_{i_s}(W)-\nabla f(W)}\\
    &= \sum_{s=0}^t\brk{\nabla f_{i_s}(W)-\nabla f(W)}-(1-\beta)\sum_{s=0}^{t-1}\sum_{j=s}^{t-1}\beta^{t-1-j}
        \brk{\nabla f_{i_s}(W)-\nabla f(W)}\\
    &= \sum_{s=0}^t\brk{\nabla f_{i_s}(W)-\nabla f(W)}-(1-\beta)\sum_{j=0}^{t-1}\beta^{t-1-j}
        \crk{\sum_{s=0}^j\brk{\nabla f_{i_s}(W)-\nabla f(W)}},
\end{aligned}
\end{equation}
where we interchanged the order of the finite sums over $0\leq s\leq j\leq t-1$ in the last equality. It then follows from \eqref{eq:sumj} that
\begin{equation}\label{eq:ema-cancel}
\begin{aligned}
    \normf{\sum_{s=0}^t\beta^{t-s}
        \prt{\nabla f_{i_s}(W)-\nabla f(W)}}
    &\leq \brk{1 + (1-\beta)\frac{1-\beta^t}{1-\beta}} \frac{1}{2}\sqrt{N\sum_{i=1}^N\normf{\nabla f_i(W)-\nabla f(W)}^2}\\
    &\leq \sqrt{N\sum_{i=1}^N\normf{\nabla f_i(W)-\nabla f(W)}^2}, \quad 0<\beta<1.
\end{aligned}
\end{equation}
For $\beta=0$, the left-hand side reduces to
$\normif{\nabla f_{i_t}(W)-\nabla f(W)}$, and the same
bound follows from \cref{lem:rr} with $A=\crki{i_t}$.

Finally, applying \eqref{eq:relative-variance} with $W = W_t$ gives
\begin{align*}
    &\sqrt{\frac{1}{N}\sum_{i=1}^N
        \normf{\nabla f_i(W_t)-\nabla f(W_t)}^2}
    \leq \sqrt{\rho\normf{\nabla f(W_t)}^2+D_0}\\
    &\leq \sqrt{\rho}\normf{\nabla f(W_t)}+\sqrt{D_0}\leq \sqrt{\rho}\norm{\nabla f(W_t)}_*+\sqrt{D_0}.
\end{align*}
Combining this estimate with \eqref{eq:ema-cancel} at $W_t$ and
the definitions of $\Gamma_1,\Gamma_0$ in \eqref{eq:heterogeneity-constants} proves
\eqref{eq:rr-descent-bound}.
\end{proof}

\begin{lemma}\label[lemma]{lem:rr}
For any $W\in\R^{m\times n}$ and any subset
$A\subseteq\crki{1,\ldots,N}$,
\begin{equation}\label{eq:subset-cancel}
    \normf{\sum_{i\in A}
        \prt{\nabla f_i(W)-\nabla f(W)}}
    \leq \frac{1}{2}\sqrt{N\sum_{i=1}^N
        \normf{\nabla f_i(W)-\nabla f(W)}^2}.
\end{equation}
\end{lemma}
\begin{proof}
Let $\1_{\crki{i\in A}} = 1$ if $i\in A$ and $0$ otherwise. Note that
\begin{align*}
    \sum_{i=1}^N\prt{\1_{\crki{i\in A}}-\frac{|A|}{N}}^2
    & = \frac{|A|(N-|A|)}{N}\leq\frac{N}{4}.
\end{align*}
It follows from $\sum_{i=1}^N(\nabla f_i(W)-\nabla f(W))=\mathbf{0}$ that
\begin{align*}
    \normf{\sum_{i\in A}
        \prt{\nabla f_i(W)-\nabla f(W)}}
    &= \normf{\sum_{i=1}^N
        \prt{\1_{\crki{i\in A}}-\frac{|A|}{N}}
        \prt{\nabla f_i(W)-\nabla f(W)}}
        \\
    &\leq \sum_{i=1}^N
        \left|\1_{\crki{i\in A}}-\frac{|A|}{N}\right|
        \normf{\nabla f_i(W)-\nabla f(W)}\\
    &\leq \sqrt{\sum_{i=1}^N
        \prt{\1_{\crki{i\in A}}-\frac{|A|}{N}}^2}
        \sqrt{\sum_{i=1}^N
            \normf{\nabla f_i(W)-\nabla f(W)}^2}\\
    &\leq \half\sqrt{N\sum_{i=1}^N
        \normf{\nabla f_i(W)-\nabla f(W)}^2}.
\end{align*}
The first inequality is the triangle inequality, and the second is the Cauchy--Schwarz inequality.
\end{proof}

\begin{lemma}\label[lemma]{lem:drift}
Let \Cref{ass:smooth} hold and $\beta\in[0,1)$. Then, the iterates $\crki{W_t}_{t\geq 0}$ generated by \Cref{alg:rr-muon} satisfy
\begin{equation}\label{eq:drift-pointwise}
\begin{aligned}
    &\frac{1-\beta}{1-\beta^{t+1}}\sum_{s=0}^t\beta^{t-s}
        \normf{\nabla f_{i_s}(W_s)-\nabla f_{i_s}(W_t)}\leq L\sqrt{r}\sum_{j=0}^{t-1}\beta^{t-j}\eta_j.
\end{aligned}
\end{equation}
\end{lemma}
\begin{proof}
Since $W_t-W_s=-\sum_{j=s}^{t-1}\eta_j\Pol(M_j)$, smoothness and \eqref{eq:polar-properties} give
\begin{align*}
    &\frac{1-\beta}{1-\beta^{t+1}}\sum_{s=0}^t\beta^{t-s}
        \normf{\nabla f_{i_s}(W_s)-\nabla f_{i_s}(W_t)}\leq\frac{L\sqrt{r}(1-\beta)}{1-\beta^{t+1}}
        \sum_{s=0}^t\beta^{t-s}\sum_{j=s}^{t-1}\eta_j\\
    &=\frac{L\sqrt{r}(1-\beta)}{1-\beta^{t+1}}
        \sum_{j=0}^{t-1}\eta_j\sum_{s=0}^j\beta^{t-s}=L\sqrt{r}\sum_{j=0}^{t-1}\eta_j\beta^{t-j}
        \frac{1-\beta^{j+1}}{1-\beta^{t+1}}\leq L\sqrt{r}\sum_{j=0}^{t-1}\eta_j\beta^{t-j}.
\end{align*}
The last inequality follows from $1-\beta^{j+1}\leq1-\beta^{t+1}$ for $j<t$. When $\beta=0$, \eqref{eq:drift-pointwise} holds because $\normif{\nabla f_{i_t}(W_t) - \nabla f_{i_t}(W_t)} = 0$.
\end{proof}

We are now ready to prove \cref{thm:recursion} by combining the previous lemmas.

\begin{proof}[Proof of \cref{thm:recursion}]
Substituting \eqref{eq:rr-descent-bound} and \eqref{eq:drift-pointwise} into \eqref{eq:descent-decomp} yields
\begin{equation}
    \label{eq:descent_prove}
    \begin{aligned}
        f(W_{t+1}) &\leq f(W_t)
        -\eta_t\prt{1-\frac{2(1-\beta)\Gamma_1}{1-\beta^{t+1}}}
            \norm{\nabla f(W_t)}_*
        +\frac{2(1-\beta)\Gamma_0\eta_t}{1-\beta^{t+1}}\\
    &\quad +2Lr\eta_t
        \sum_{j=0}^{t-1}\beta^{t-j}\eta_j +\frac{2\sqrt{r}\eta_t\beta^{t+1}}{1-\beta^{t+1}}
        \normf{M_{\mathrm{init}}}+\frac{Lr}{2}\eta_t^2,\quad t\geq 0,
    \end{aligned}
\end{equation}
which is \eqref{eq:descent}. Using $1/(1-\beta^{t+1})=1+[\beta^{t+1}/(1-\beta^{t+1})]$, we have 
\begin{equation}
    \label{eq:separate}
    \frac{2(1-\beta)\Gamma_0\eta_t}{1-\beta^{t+1}} = 2(1-\beta)\Gamma_0\eta_t\prt{1+\frac{\beta^{t+1}}{1-\beta^{t+1}}},\quad \frac{2(1-\beta)\Gamma_1\eta_t}{1-\beta^{t+1}} = 2(1-\beta)\Gamma_1\eta_t\prt{1+\frac{\beta^{t+1}}{1-\beta^{t+1}}}.
\end{equation}

To show \eqref{eq:weighted-finite}, assume $2(1-\beta)\Gamma_1<1$.
Fix $T\geq1$. Summing \eqref{eq:descent_prove} over $t=0,\ldots,T-1$ and invoking \eqref{eq:separate} yields
\begin{equation}\label{eq:unified-descent}
\begin{aligned}
    &\brk{1-2(1-\beta)\Gamma_1}
        \sum_{t=0}^{T-1}\eta_t\norm{\nabla f(W_t)}_*\\
    &\leq f(W_0)-f(W_T)+2(1-\beta)\Gamma_0
        \sum_{t=0}^{T-1}\eta_t+\frac{Lr}{2}\sum_{t=0}^{T-1}\eta_t^2
        +2Lr\sum_{t=0}^{T-1}\eta_t\sum_{j=0}^{t-1}\beta^{t-j}\eta_j\\
    &\quad+2\sqrt{r}\normf{M_{\mathrm{init}}}
        \sum_{t=0}^{T-1}\frac{\eta_t\beta^{t+1}}{1-\beta^{t+1}} +2(1-\beta)\sum_{t=0}^{T-1}
        \frac{\eta_t\beta^{t+1}}{1-\beta^{t+1}}
        \prt{\Gamma_1\norm{\nabla f(W_t)}_*+\Gamma_0}.
\end{aligned}
\end{equation}

We next estimate the accumulated terms.
First, using $\eta_t\leq\eta_j$ for $j<t$ and interchanging the sums, we obtain
\begin{equation}\label{eq:drift-sum}
\begin{aligned}
    &\sum_{t=0}^{T-1}\eta_t
        \sum_{j=0}^{t-1}\beta^{t-j}\eta_j\leq \sum_{j=0}^{T-1}\eta_j^2
        \sum_{t=j+1}^{T-1}\beta^{t-j}=\sum_{j=0}^{T-1}\eta_j^2
        \frac{\beta(1-\beta^{T-1-j})}{1-\beta} \leq \frac{\beta}{1-\beta}\sum_{j=0}^{T-1}\eta_j^2.
\end{aligned}
\end{equation}
Next, $1-\beta^{t+1}\geq1-\beta$ and $\eta_t\leq\eta_0$ give
\begin{equation}\label{eq:init-sum}
\begin{aligned}
    \sum_{t=0}^{T-1}\frac{\eta_t\beta^{t+1}}{1-\beta^{t+1}}
    &\leq \frac{\eta_0}{1-\beta}\sum_{t=0}^{T-1}\beta^{t+1}\leq \frac{\eta_0\beta}{(1-\beta)^2}.
\end{aligned}
\end{equation}
Consequently, 
\begin{equation}\label{eq:init-smooth-sum}
    \sum_{t=0}^{T-1}\prt{\frac{2\sqrt{r}\eta_t\beta^{t+1}}{1-\beta^{t+1}}\normf{M_{\mathrm{init}}} + \frac{Lr}{2}\eta_t^2 + 2Lr\eta_t\sum_{j=0}^{t-1}\beta^{t-j}\eta_j}
    \leq \frac{2\sqrt{r}\eta_0\beta}{(1-\beta)^2}\normf{M_{\mathrm{init}}}+\prt{\frac{1}{2} + \frac{2\beta}{1-\beta}}Lr\sum_{t=0}^{T-1}\eta_t^2.
\end{equation}

For the last sum in \eqref{eq:unified-descent}, smoothness and \eqref{eq:step-length} give
\begin{equation}\label{eq:gradient-growth-time}
\begin{aligned}
    &\left|\norm{\nabla f(W_{t+1})}_*-\norm{\nabla f(W_t)}_*\right|
    \leq \norm{\nabla f(W_{t + 1}) - \nabla f(W_t)}_*\leq \sqrt{r}\normf{\nabla f(W_{t + 1}) - \nabla f(W_t)}\leq Lr\eta_t,\\
    &\norm{\nabla f(W_t)}_*
    \leq\norm{\nabla f(W_0)}_*+Lr\eta_0t.
\end{aligned}
\end{equation}
It follows from $\eta_t\leq\eta_0$ and $1-\beta^{t+1}\geq1-\beta$ that 
\begin{equation}\label{eq:transient-sum}
\begin{aligned}
    2(1-\beta)\sum_{t=0}^{T-1}
        \frac{\eta_t\beta^{t+1}}{1-\beta^{t+1}}
        \prt{\Gamma_1\norm{\nabla f(W_t)}_*+\Gamma_0}&\leq2\eta_0\sum_{t=0}^{\infty}\beta^{t+1}
        \prt{\Gamma_1\norm{\nabla f(W_0)}_*
            +\Gamma_1 Lr\eta_0t+\Gamma_0}\\
    &=\frac{2\eta_0\beta}{1-\beta}
        \prt{\Gamma_1\norm{\nabla f(W_0)}_*+\Gamma_0}
        +\frac{2\Gamma_1 Lr\eta_0^2\beta^2}{(1-\beta)^2},
\end{aligned}
\end{equation}
where we used
$\sum_{t=0}^{\infty}\beta^{t+1}=\beta/(1-\beta)$ and
$\sum_{t=0}^{\infty}t\beta^{t+1}=\beta^2/(1-\beta)^2$.
These estimates also cover $\beta=0$. Substituting
\eqref{eq:init-smooth-sum} and \eqref{eq:transient-sum} into
\eqref{eq:unified-descent} and using
$f(W_0)-f(W_T)\leq\Delta_f$ yields
\begin{equation}\label{eq:unified-sum}
\begin{aligned}
    &\prt{1-2(1-\beta)\Gamma_1}
        \sum_{t=0}^{T-1}\eta_t\norm{\nabla f(W_t)}_*\leq\Delta_f+Lr\prt{\frac{1}{2}+\frac{2\beta}{1-\beta}}
        \sum_{t=0}^{T-1}\eta_t^2
        +2(1-\beta)\Gamma_0\sum_{t=0}^{T-1}\eta_t\\
    &\quad+\frac{2\eta_0\beta}{(1-\beta)^2}
        \brk{\sqrt{r}\normf{M_{\mathrm{init}}}
            +(1-\beta)\prt{\Gamma_1\norm{\nabla f(W_0)}_*+\Gamma_0}}+\frac{2\Gamma_1 Lr\eta_0^2\beta^2}{(1-\beta)^2}.
\end{aligned}
\end{equation}
Since $2(1-\beta)\Gamma_1<1$, dividing by $(1-2(1-\beta)\Gamma_1)\sum_{t=0}^{T-1}\eta_t$ proves \eqref{eq:weighted-finite}.

\end{proof}

\subsection{Proof of Corollary \ref{cor:complexity}}

With $M_{\mathrm{init}}=\mathbf{0}$, a constant stepsize $\eta_t=\eta$, and $2(1-\beta)\Gamma_1\leq 1/2$, \eqref{eq:weighted-finite} implies
\begin{equation}\label{eq:constant-step-finite}
\begin{aligned}
    \frac{1}{T}\sum_{t=0}^{T-1}\norm{\nabla f(W_t)}_*
    &\leq\frac{2\Delta_f}{T\eta}
        +\frac{4\beta\prt{\Gamma_1\norm{\nabla f(W_0)}_*+\Gamma_0}}{(1-\beta)T} +\frac{4\Gamma_1 Lr\eta\beta^2}{(1-\beta)^2T}
        +\frac{4Lr\eta}{(1-\beta)} +4(1-\beta)\Gamma_0.
\end{aligned}
\end{equation}

The choice of parameters in \eqref{eq:horizon-parameters} leads to
\begin{equation}
    \label{eq:stepsize-in}
    \begin{aligned}
        \frac{2\Delta_f}{T\eta}&=\frac{4Lr\eta}{1-\beta}=2\sqrt{\frac{2Lr\Delta_f}{(1-\beta)T}},\\
        \frac{1}{1-\beta} &= 1+4\Gamma_1+\prt{\frac{2\Gamma_0^2T}{Lr\Delta_f}}^{1/3},\quad 1-\beta\leq \prt{\frac{Lr\Delta_f}{2\Gamma_0^2T}}^{1/3},\\
        \frac{1}{(1-\beta)^{3/2}} &= \prt{1+4\Gamma_1+\prt{\frac{2\Gamma_0^2T}{Lr\Delta_f}}^{1/3}}^{3/2}\leq 2^{3/2-1}\brk{(1 + 4\Gamma_1)^{3/2} + \sqrt{\frac{2\Gamma_0^2T}{Lr\Delta_f}}},
    \end{aligned}
\end{equation}
where we used the convexity of $x\mapsto x^{3/2}$ on $[0,\infty)$. Substituting \eqref{eq:stepsize-in} into \eqref{eq:constant-step-finite} yields 
\begin{equation}\label{eq:horizon-rate}
\begin{aligned}
    \frac{1}{T}\sum_{t=0}^{T-1}\norm{\nabla f(W_t)}_*
    &\leq10\prt{\frac{Lr\Gamma_0\Delta_f}{T}}^{1/3}+4\sqrt{\frac{2(1+4\Gamma_1)Lr\Delta_f}{T}}
        +\frac{4\prt{2\Gamma_0^2/(Lr\Delta_f)}^{1/3}\prt{\Gamma_1\norm{\nabla f(W_0)}_*+\Gamma_0}}{T^{2/3}}\\
    &\quad+\frac{6(1+4\Gamma_1)\prt{\Gamma_1\norm{\nabla f(W_0)}_*+\Gamma_0}}{T}
        +\frac{4\Gamma_1(1+4\Gamma_1)^{3/2}\sqrt{Lr\Delta_f}}{T^{3/2}},
\end{aligned}
\end{equation}
which gives \eqref{eq:horizon-rate-order}. Requiring each of the five terms in \eqref{eq:horizon-rate} to be at most $\varepsilon/5$ gives \eqref{eq:complexity-horizon}.

\subsection{Proof of Corollary \ref{cor:exact}}
\label[appendix]{app:strong-growth}

We first prove $f(W_t)\to f_\infty$. Since $D_0 = 0$ and $M_{\mathrm{init}} = \mathbf{0}$, \eqref{eq:descent} becomes 
\begin{equation}
    \label{eq:descent_Gamma0}
    \begin{aligned}
        f(W_{t+1}) &\leq f(W_t)
        -\eta_t\prt{1-\frac{2(1-\beta)\Gamma_1}{1-\beta^{t+1}}}
            \norm{\nabla f(W_t)}_*
        +2Lr\eta_t
        \sum_{j=0}^{t-1}\beta^{t-j}\eta_j +\frac{Lr}{2}\eta_t^2,\quad t\geq 0.
    \end{aligned}
\end{equation}
Because $2(1-\beta)\Gamma_1<1$ and $\beta^{t+1}\to0$, there is an integer
$T_0\geq0$ such that $1-2(1-\beta)\Gamma_1/(1-\beta^{t+1})>0$ for all
$t\geq T_0$. Hence
\begin{equation}\label{eq:quasi-monotone}
    f(W_{t+1})\leq f(W_t)+a_t,\quad
    a_t:=2Lr\eta_t\sum_{j=0}^{t-1}\beta^{t-j}\eta_j+\frac{Lr}{2}\eta_t^2\geq0,\qquad t\geq T_0.
\end{equation}
Letting $T\to\infty$ in \eqref{eq:drift-sum} and using
$\sum_t\eta_t^2<\infty$ gives $\sum_t a_t<\infty$. Hence, $\crki{f(W_t)-\sum_{s=T_0}^{t-1}a_s}_t$ is nonincreasing for $t\geq T_0$ and bounded below by $f^*-\sum_{s\geq T_0}a_s$, so it converges. Since the partial
sums of $a_s$ also converge, $f(W_t)$ converges to a finite limit
$f_\infty\geq f^*$.

We now show that $\normi{\nabla f(W_t)}_*\to0$. The estimate
\eqref{eq:unified-sum} and $\sum_t\eta_t^2<\infty$ yield
\begin{equation}\label{eq:weighted-summable}
\begin{aligned}
    \sum_{t=0}^{\infty}\eta_t\norm{\nabla f(W_t)}_*
    &\leq\frac{\Delta_f+Lr
        \prt{\frac{1}{2}+\frac{2\beta}{1-\beta}}
        \sum_{t=0}^{\infty}\eta_t^2}{1-2(1-\beta)\Gamma_1}\\
    &\quad+\frac{2\eta_0\beta\Gamma_1\norm{\nabla f(W_0)}_*}
        {(1-\beta)\prt{1-2(1-\beta)\Gamma_1}}
        +\frac{2\Gamma_1 Lr\eta_0^2\beta^2}
        {(1-\beta)^2\prt{1-2(1-\beta)\Gamma_1}}
     <\infty.
\end{aligned}
\end{equation}

Since $\sum_t\eta_t=\infty$, \eqref{eq:weighted-summable} implies that
$\liminf_t\norm{\nabla f(W_t)}_*=0$. We prove $\normi{\nabla f(W_t)}_*\to0$ by contradiction. If $\normi{\nabla f(W_t)}_*\not\to0$, then there exists an $\varepsilon>0$ and infinitely many disjoint index intervals $[u,v)$ such that
\[
    \norm{\nabla f(W_u)}_*\geq2\varepsilon,\quad
    \norm{\nabla f(W_v)}_*\leq\varepsilon,\quad
    \norm{\nabla f(W_t)}_*\geq\varepsilon\quad(u\leq t<v).
\]
By \eqref{eq:gradient-growth-time}, every such interval satisfies
\[
    \sum_{t=u}^{v-1}\eta_t\norm{\nabla f(W_t)}_*
    \geq\varepsilon\sum_{t=u}^{v-1}\eta_t
    \geq\frac{\varepsilon}{Lr}
        \prt{\norm{\nabla f(W_u)}_*-\norm{\nabla f(W_v)}_*}
    \geq\frac{\varepsilon^2}{Lr}.
\]
Summing over the disjoint intervals contradicts
\eqref{eq:weighted-summable}. Therefore
$\norm{\nabla f(W_t)}_*\to0$. By the continuity of $\nabla f(W)$, every accumulation point is stationary. 
%\end{proof}

\section{Proof of the Failure Result}
\label[appendix]{app:smooth-counterexample-proof}

\textbf{Part (i).}
Every $g_i$ is nonnegative, convex, and continuously differentiable with an $N$-Lipschitz derivative. Since
$h'(u)=\max\crki{-1,\min\crki{u,1}}$, we obtain
\begin{equation}\label{eq:g_gradient}
    g'(x)=h'(x+2)+\frac{N-1}{N}h'(x-2)
    = \begin{cases}
        -\frac{2N-1}{N}, & x<-3,\\[4pt]
        x+1+\frac{1}{N}, & -3\leq x<-1,\\[4pt]
        \frac{1}{N}, & -1\leq x\leq 1,\\[4pt]
        1+\frac{N-1}{N}(x-2), & 1<x\leq 3,\\[4pt]
        \frac{2N-1}{N}, & x>3.
    \end{cases}
\end{equation}
Hence, $g'$ is nondecreasing, and $g'(x)=0$ holds only at
\begin{equation}\label{eq:huber-minimizer}
    x_*=-1-\frac{1}{N}.
\end{equation}
Since $g'<0$ on $(-\infty,x_*)$ and $g'>0$ on $(x_*,\infty)$, the point $x_*$ is the unique minimizer of $g$.
Moreover, $g_2=\cdots=g_N$, so 
$g_1'(x)-g'(x)=\frac{N-1}{N}[g_1'(x)-g_2'(x)]$, and
$g_i'(x)-g'(x)=-\frac{1}{N}[g_1'(x)-g_2'(x)]$ for $i\geq2$. Hence
\begin{equation}\label{eq:huber-growth}
\begin{aligned}
    \frac{1}{N}\sum_{i=1}^{N}|g_i'(x)|^2-|g'(x)|^2
    &=\frac{1}{N}\sum_{i=1}^{N}|g_i'(x)-g'(x)|^2
     =\frac{1}{N}\prt{\frac{(N-1)^2}{N^2}+\frac{N-1}{N^2}}|g_1'(x)-g_2'(x)|^2\\
    &=\frac{N-1}{N^2}|g_1'(x)-g_2'(x)|^2\leq\frac{(N-1)(N+1)^2}{N^2},\quad\forall x.
\end{aligned}
\end{equation}
Here, the last inequality holds because $|g_1'(x)|\leq N$ and $|g_2'(x)|\leq 1$.
We now turn to $f_i$ and $f$. Each $f_i$ is nonnegative and convex as a sum of two nonnegative convex functions, and
coercive because $g_i$ is coercive and the quadratic term is coercive in the
other entries.
Let $E_{11}:=e_1e_1^\T$ with $e_1:= (1,0)^\T$. By \eqref{eq:huber-matrix},
\begin{equation}\label{eq:huber-matrix-gradient}
    \nabla f_i(W)=W+\brk{g_i'([W]_{11})-[W]_{11}}E_{11},\quad
    \nabla f(W)=W+\brk{g'([W]_{11})-[W]_{11}}E_{11}.
\end{equation}
For $W,Z\in\R^{2\times2}$,
\begin{equation}
    \label{eq:part1_smooth}
     \normf{\nabla f_i(W)-\nabla f_i(Z)}^2 =|g_i'([W]_{11})-g_i'([Z]_{11})|^2
     +\sum_{(a,b)\neq(1,1)}\prt{[W]_{ab}-[Z]_{ab}}^2 \leq N^2\normf{W-Z}^2,
\end{equation}
where we used that $g_i'$ is $N$-Lipschitz and $N\geq1$. Since $f(W)\geq g([W]_{11})$, we have
$f^*\geq g(x_*)>-\infty$, so \Cref{ass:smooth} holds with $L=N$. By
\eqref{eq:huber-matrix-gradient}, the entries $[\nabla f_i(W)]_{ab}=[W]_{ab}$
with $(a,b)\neq(1,1)$ are the same for all $i$ and cancel in $\nabla f_i(W)-\nabla f(W)=\prt{g_i'([W]_{11})-g'([W]_{11})}E_{11}$. The
identity in \eqref{eq:relative-variance} and \eqref{eq:huber-growth} then
give
\begin{equation}
    \label{eq:part1_growth}
    \frac{1}{N}\sum_{i=1}^{N}\normf{\nabla f_i(W)}^2-\normf{\nabla f(W)}^2
    =\frac{1}{N}\sum_{i=1}^{N}|g_i'([W]_{11})-g'([W]_{11})|^2
    \leq\frac{(N-1)(N+1)^2}{N^2},
\end{equation}
so \Cref{ass:growth} holds with $D_1=1$ and $D_0=(N-1)(N+1)^2/N^2>0$.

\textbf{Part (ii).}
We prove \eqref{eq:huber-failure}. We first show that the evolution of $W_t$ under \cref{alg:rr-muon} reduces to the scalar sequence $x_t$ defined in \eqref{eq:scalar-muon}.
For $W=xE_{11}$, \eqref{eq:huber-matrix-gradient} gives
$\nabla f_i(W)=g_i'(x)E_{11}$ and $\nabla f(W)=g'(x)E_{11}$. Let $(x_t,m_t)$ be generated by
\begin{equation}\label{eq:scalar-muon}
    m_t=\beta m_{t-1}+(1-\beta)g_{i_t}'(x_t),\qquad
    x_{t+1}=x_t-\eta_t\sgn(m_t),\qquad x_0=m_{-1}=0,
\end{equation}
with $\sgn(0)=0$. We show by induction that $W_t=x_tE_{11}$ and
$M_{t-1}=m_{t-1}E_{11}$ for all $t\geq0$. The case $t=0$ holds because $W_0=M_{-1}=0$. If the claim holds for $t$, then $\nabla f_{i_t}(W_t)=g_{i_t}'(x_t)E_{11}$, and the momentum update of \cref{alg:rr-muon} together with \eqref{eq:scalar-muon} gives
\[
    M_t=\beta m_{t-1}E_{11}+(1-\beta)g_{i_t}'(x_t)E_{11}=\brk{\beta m_{t-1}+(1-\beta)g_{i_t}'(x_t)}E_{11} = m_t E_{11}.
\]
If $m_t\neq0$, then $M_t$ has rank one and admits the compact singular value
decomposition $M_t=U_1\Sigma_1V_1^\T$ with
\[
    U_1=\sgn(m_t)\,e_1\in\R^{2\times1},\qquad
    \Sigma_1=|m_t|\in\R^{1\times1},\qquad
    V_1=e_1\in\R^{2\times1},
\]
so \eqref{eq:polar} gives $\Pol(M_t)=U_1V_1^\T=\sgn(m_t)E_{11}$. This also holds for $m_t=0$ because $\Pol(\mathbf{0})=\mathbf{0}$. Hence
\begin{align*}
    W_{t+1}=W_t-\eta_t\Pol(M_t)=\brk{x_t - \eta_t\sgn(m_t)}E_{11} = x_{t+1}E_{11},
\end{align*}
which completes the induction.
Consequently, $\nabla f(W_t)=g'(x_t)E_{11}$ and $\norm{\nabla f(W_t)}_*=|g'(x_t)|$.
Moreover, $f^*\geq g(x_*)$ by Part (i) and $f(x_*E_{11})=g(x_*)$ by \eqref{eq:huber-matrix}, so $f^*=g(x_*)$ and $f(W_t)-f^*=g(x_t)-g(x_*)$.

To prove \eqref{eq:huber-failure}, it therefore suffices to show $|g'(x_t)|\geq1/N$ and
$g(x_t)-g(x_*)\geq(N+1)/(2N^2)$ for all $t$. We prove this through the
chain
\begin{equation}\label{eq:huber-chain}
    x_{kN}\geq0\ \ \forall k
    \quad\Longrightarrow\quad
    x_t\geq-\frac{1}{2}\ \ \forall t
    \quad\Longrightarrow\quad
    |g'(x_t)|\geq\frac{1}{N}\ \text{ and }\ g(x_t)-g(x_*)\geq\frac{N+1}{2N^2}\ \ \forall t.
\end{equation}

In the following, Steps 1 and 2 prove the two implications given $x_{kN}\geq0$ for all $k\geq 0$, which is established in Step 3. 

\emph{Step 1: $x_{kN}\geq0$ for all $k$ implies $x_t\geq-1/2$ for all $t$.} Since $x_{kN+j} = x_{kN} - \sum_{i=0}^{j-1} \eta_{kN+i} \sgn(m_{kN+i})$ and $\eta_{kN + i}=\eta_{kN}\leq 1/(2N)$ for $0\leq i<N$, we have
\begin{equation}\label{eq:huber-epoch-displacement}
    |x_{kN+j}-x_{kN}|\leq N\eta_{kN}\leq\frac{1}{2},\qquad 0\leq j\leq N.
\end{equation}
Hence, if $x_{kN}\geq0$ for all $k$, then $x_{kN+j}\geq x_{kN}-1/2\geq -1/2$. Since every $t$ can be written as $t=kN+j$ with $0\leq j<N$, we obtain $x_t\geq-1/2$ for all $t$.

\emph{Step 2: $x_t\geq-1/2$ for all $t$ implies the two bounds.} By \eqref{eq:g_gradient}, $g'$ is
nondecreasing, so $x_t\geq-1/2$ gives $g'(x_t)\geq g'(-1/2)=1/N$. Since $g'\geq0$ on $[x_*,\infty)$ and $x_*=-1-1/N<-1/2$, the function $g$ is nondecreasing on $[-1/2,\infty)$. Consequently, $g(x_t)\geq g(-1/2)$. 
Moreover, we have
\[
\begin{aligned}
    g\prt{-\frac{1}{2}}-g(x_*)
    &=\brk{h\prt{\frac{3}{2}}+\frac{N-1}{N}h\prt{-\frac{5}{2}}}
     -\brk{h\prt{1-\frac{1}{N}}+\frac{N-1}{N}h\prt{-3-\frac{1}{N}}}\\
    &=\brk{1+\frac{2(N-1)}{N}}
     -\brk{\frac{(N-1)^2}{2N^2}+\frac{N-1}{N}\prt{\frac{5}{2}+\frac{1}{N}}}
     =\frac{N+1}{2N^2}.
\end{aligned}
\]
Hence, $g(x_t)-g(x_*)\geq g(-1/2)-g(x_*)=(N+1)/(2N^2)$.

\emph{Step 3: $x_{kN}\geq0$ for all $k$.} We prove the claim by induction on $k$.
The base case is $x_0=0$. Suppose $x_{kN}\geq0$. Then \eqref{eq:huber-epoch-displacement} yields
\[
    -\frac{1}{2}\leq x_{kN} - \frac{1}{2} \leq x_{kN+j}\leq x_{kN}+ \frac{1}{2},\qquad 0\leq j\leq N.
\]
If $x_{kN}\geq1/2$, then $x_{(k+1)N}\geq0$. The induction step is complete in this case. 
Otherwise, $0\leq x_{kN}<1/2$, and then $x_{kN+j}\in[-1/2,1]\subset[-1,1]$ for $0\leq j<N$. In this case, it suffices to show that epoch $k$ increases $x_{kN}$, that is, $x_{(k+1)N}>x_{kN}$. We use three facts.

(a) \emph{Inputs on $[-1,1]$.} For $x\in[-1,1]$, we have 
\[
    g_1'(x)=Nh'(x+2)=N,\qquad g_i'(x)=h'(x-2)=-1\quad(2\leq i\leq N).
\]
Hence, in an epoch whose evaluation
points all lie in $[-1,1]$, the input $g_{i_t}'(x_t)$ equals $N$ when $i_t=1$ and $-1$ otherwise. Since each epoch visits every component exactly once, such an epoch contains exactly one large input $g_1'=N$ and $N-1=4K$ inputs $-1$, for every sequence of permutations.

(b) \emph{Momentum.} The global gradient bounds $|g_1'(x)|\leq N$ and
$|g_i'(x)|\leq1$ for $i\geq2$ imply $m_t\in[-N,N]$ for every $t$ since $m_t = (1-\beta)\sum_{s=0}^{t} \beta^{t-s} g_{i_s}'(x_s)$. If the
inputs at times $s=t+1,\ldots,t+K$ all equal $-1$, then for every $m_t\in[-N,N]$,
\begin{equation}\label{eq:huber-forgetting}
    m_{t+K}=\beta^K m_t-(1-\beta)\sum_{j=0}^{K-1}\beta^j=-1+\beta^K(m_t+1)\leq-1+(N+1)\beta^K<0,
\end{equation}
where the last inequality uses $(N+1)\beta^K=(4K+2)\beta^K<1$. The momentum stays negative as long as the subsequent inputs equal $-1$,
since $\beta m-(1-\beta)<0$ whenever $m<0$. The large input can make it nonnegative again. By (a), this happens at most once per epoch, and (c)
accounts for it.

\begin{figure}[!htbp]
    \centering
    \begin{tikzpicture}[x=0.95cm,y=1cm,
        cell/.style={draw, minimum width=0.95cm, minimum height=0.6cm,
            inner sep=0pt, font=\small},
        maybe/.style={cell, fill=gray!25}]
        % Case a >= K: a = 3.
        \foreach \j/\v/\s in {0/-1/maybe,1/-1/cell,2/-1/cell,3/N/maybe,
                4/-1/maybe,5/-1/cell,6/-1/cell,7/-1/cell,8/-1/cell}{
            \node[\s] at (\j,0) {$\v$};
            \node[font=\scriptsize, gray] at (\j,-0.5) {$\j$};
        }
        \draw[decorate,decoration={brace,amplitude=4pt}] (-0.45,0.38)
            -- node[above=5pt,font=\small]{$a=3$} (2.45,0.38);
        \draw[decorate,decoration={brace,amplitude=4pt}] (3.55,0.38)
            -- node[above=5pt,font=\small]{$4K-a=5$} (8.45,0.38);
        \node[anchor=east,font=\small] at (-0.7,0) {$a\geq K$};
        % Case a < K: a = 1.
        \begin{scope}[yshift=-2.1cm]
        \foreach \j/\v/\s in {0/-1/maybe,1/N/maybe,2/-1/maybe,3/-1/cell,
                4/-1/cell,5/-1/cell,6/-1/cell,7/-1/cell,8/-1/cell}{
            \node[\s] at (\j,0) {$\v$};
            \node[font=\scriptsize, gray] at (\j,-0.5) {$\j$};
        }
        \draw[decorate,decoration={brace,amplitude=4pt}] (-0.45,0.38)
            -- node[above=5pt,font=\small]{$a=1$} (0.45,0.38);
        \draw[decorate,decoration={brace,amplitude=4pt}] (1.55,0.38)
            -- node[above=5pt,font=\small]{$4K-a=7$} (8.45,0.38);
        \node[anchor=east,font=\small] at (-0.7,0) {$a<K$};
        \end{scope}
        % Case 4K - a < K: a = 7.
        \begin{scope}[yshift=-4.2cm]
        \foreach \j/\v/\s in {0/-1/maybe,1/-1/cell,2/-1/cell,3/-1/cell,
                4/-1/cell,5/-1/cell,6/-1/cell,7/N/maybe,8/-1/maybe}{
            \node[\s] at (\j,0) {$\v$};
            \node[font=\scriptsize, gray] at (\j,-0.5) {$\j$};
        }
        \draw[decorate,decoration={brace,amplitude=4pt}] (-0.45,0.38)
            -- node[above=5pt,font=\small]{$a=7$} (6.45,0.38);
        \draw[decorate,decoration={brace,amplitude=4pt}] (7.55,0.38)
            -- node[above=5pt,font=\small]{$4K-a=1$} (8.45,0.38);
        \node[anchor=east,font=\small] at (-0.7,0) {$4K-a<K$};
        \end{scope}
    \end{tikzpicture}
    \caption{One epoch whose gradient-evaluation points lie in $[-1,1]$, for
    $K=2$ ($N=9$), with $a=3$ (top), $a=1$ (middle), and $a=7$ (bottom) inputs
    $-1$ before the large input $N$. The top row has $a\geq K$, the middle row
    $a<K$, and the bottom row $4K-a<K$; the top and middle rows have
    $4K-a\geq K$. The braces mark the runs of inputs $-1$, and the small
    numbers index the updates of the epoch. Shaded updates may have
    nonnegative momentum: at most $K-1$ before the large input and at most $K$
    from it on.}
    \label{fig:huber-epoch}
\end{figure}
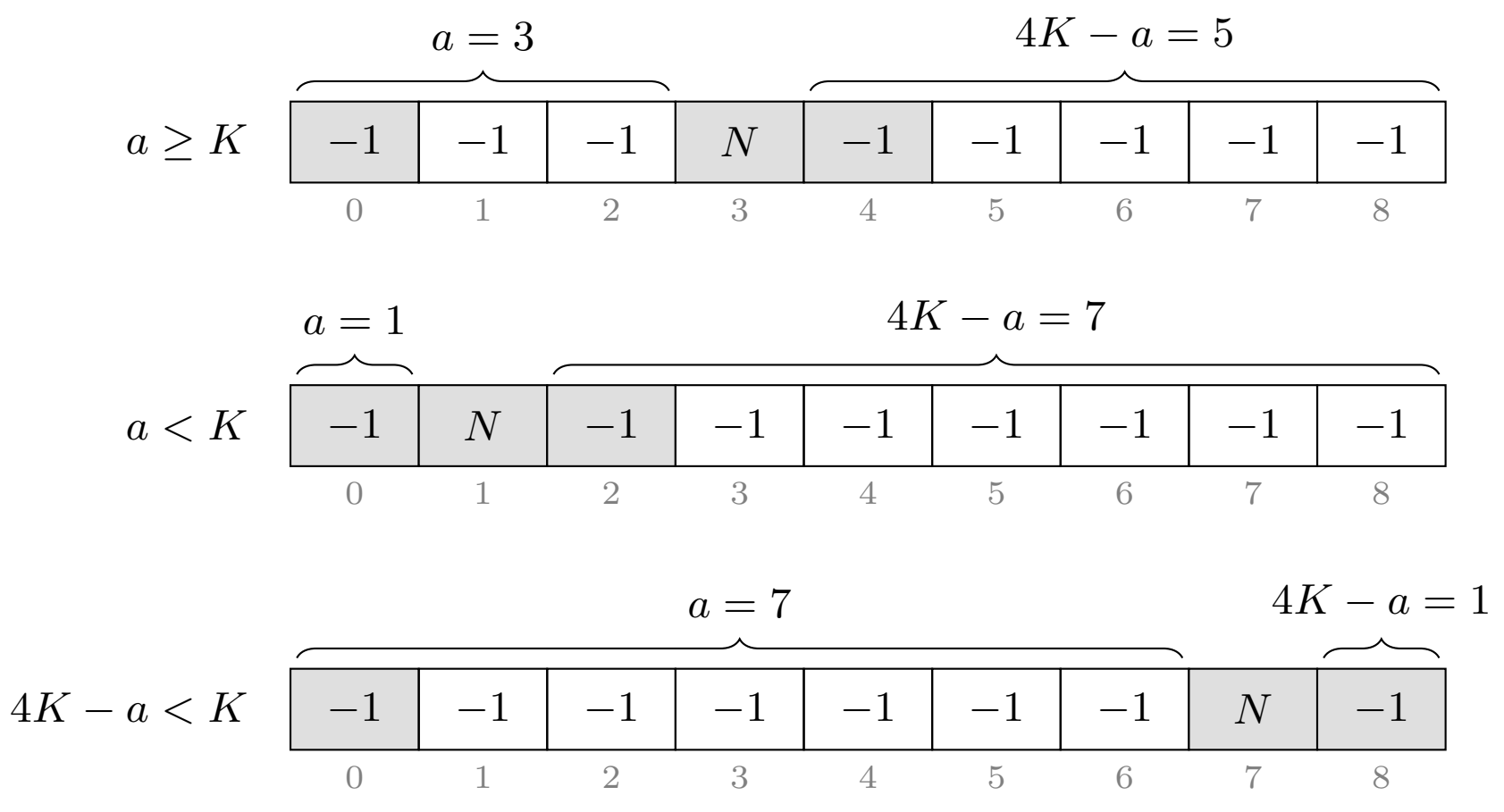

(c) \emph{Increase over the epoch.} Consider an epoch whose gradient-evaluation points all lie in $[-1,1]$. By (a), its inputs are $4K$ copies of $-1$ and one large input $g_1'=N$.
The large input splits the inputs $-1$ into a run of length $a$ before it and a run of length $4K-a$ after it (\Cref{fig:huber-epoch}). We count the updates
with nonnegative momentum in each part.
\begin{itemize}
    \item \emph{Before the large input.} If $a\geq K$
    (\Cref{fig:huber-epoch}, top), the first $K$ inputs of the epoch satisfy
    the premise of (b) with $t=kN-1$, so the momentum is negative from the
    $K$th update of the epoch until the large input. If $a<K$
    (\Cref{fig:huber-epoch}, middle), only $a\leq K-1$ updates precede the large input. In both cases, at most $K-1$ updates before the large input have nonnegative momentum.
    \item \emph{From the large input on.} Let the large input occur at time $s$. If $4K-a\geq K$ (\Cref{fig:huber-epoch}, top and middle), the inputs at times $s+1,\ldots,s+K$ satisfy the premise of (b) with $t=s$, so the
    momentum is negative from time $s+K$ until the end of the epoch. If $4K-a<K$ (\Cref{fig:huber-epoch}, bottom), only $4K-a+1\leq K$ updates remain from time $s$ on. In both cases, at most $K$ updates from time $s$
    on have nonnegative momentum.
\end{itemize}
Hence, at most $2K$ of the $N$ updates of the epoch have nonnegative
momentum, and at least $N-2K$ have negative momentum. Summing over the epoch according to \eqref{eq:scalar-muon}, we obtain
\begin{equation}\label{eq:huber-epoch-drift}
    x_{(k+1)N}-x_{kN}
    =-\eta_{kN}\sum_{j=0}^{N-1}\sgn(m_{kN+j})\geq \eta_{kN}\prt{(N-2K)-2K}=\eta_{kN}>0,
\end{equation}
where the last equality uses $N-4K=1$. This completes the induction for the case $0\leq x_{kN}<1/2$, so $x_{kN}\geq0$ for all $k$. This establishes \eqref{eq:huber-chain} and hence \eqref{eq:huber-failure}.

\textbf{Part (iii).}
Let $\tilde{f}_i:=f_i/N$ and $\tilde{f}:=f/N$. Then, $\nabla\tilde{f}_i=\nabla f_i/N$ and $\nabla\tilde{f}=\nabla f/N$. By Part (i), each $\tilde{f}_i$ is convex with a $1$-Lipschitz gradient, and $\tilde{f}^*=f^*/N>-\infty$. Applying \eqref{eq:part1_growth} to $\tilde{f}_i$ and $\tilde{f}$ yields
\[
    \frac{1}{N}\sum_{i=1}^{N}\normf{\nabla\tilde{f}_i(W)}^2
    -\normf{\nabla\tilde{f}(W)}^2 = \frac{1}{N^3}\sum_{i=1}^{N} \normf{\nabla f_i(W)}^2 - \frac{1}{N^2}\normf{\nabla f(W)}^2 \leq \frac{(N-1)(N+1)^2}{N^4},\quad\forall W\in\R^{2\times2},
\]
so \Cref{ass:growth} holds with $D_1=1$ and $D_0=(N-1)(N+1)^2/N^4\leq 1$ for $N\geq 2$. 

We next show that \cref{alg:rr-muon} generates the same iterates on $\tilde{f}_1,\ldots,\tilde{f}_N$ as on $f_1,\ldots,f_N$. Let $(W_t,M_t)$ and $(\tilde W_t,\tilde M_t)$ denote the iterates on the two families, generated from $W_0=\tilde W_0=0$ and $M_{-1}=\tilde M_{-1}=0$ with the same permutations and the same stepsizes. These stepsizes are admissible for both families because the stepsize choice in \eqref{eq:huber-stepsizes} depends only on $N$. We show by induction that $\tilde W_t=W_t$ and $\tilde M_{t-1}=M_{t-1}/N$ for all $t\geq0$. The case $t=0$ holds by the initialization. If the claim holds for $t$, then $\nabla\tilde{f}_{i_t}(\tilde W_t)=\nabla f_{i_t}(W_t)/N$, and the momentum update of \cref{alg:rr-muon} gives
\[
    \tilde M_t=\beta\frac{M_{t-1}}{N}+(1-\beta)\frac{\nabla f_{i_t}(W_t)}{N}=\frac{M_t}{N}.
\]
Since $\Pol(aM)=\Pol(M)$ for $a>0$ by \eqref{eq:polar-properties}, we have $\Pol(\tilde M_t)=\Pol(M_t)$, and hence
\[
    \tilde W_{t+1}=\tilde W_t-\eta_t\Pol(\tilde M_t)=W_t-\eta_t\Pol(M_t)=W_{t+1}.
\]
This completes the induction.

Finally, $\nabla\tilde{f}(\tilde W_t)=\nabla f(W_t)/N$ for all $t\geq0$, so \eqref{eq:huber-failure} gives $\normi{\nabla\tilde{f}(\tilde W_t)}_* = \normi{\nabla f(W_t)/N}_* \geq1/N^2$.

\section{Implementation Details}\label[appendix]{app:implementation}

This appendix records a PyTorch configuration that realizes
\cref{alg:rr-muon} up to the Newton--Schulz approximation of $\Pol$.
%,and relates the two momentum conventions used by PyTorch.

%\textbf{PyTorch configuration.}
The setup is similar in structure to Moonlight's official illustrative
training script, which combines \texttt{DataLoader(..., shuffle=True)} with a Muon optimizer. Here we use \texttt{torch.optim.Muon}
with Nesterov momentum disabled and no weight decay to match
\cref{alg:rr-muon}. For a trainable matrix tensor \texttt{W}, the
PyTorch~2.14 configuration combines RR data loading with zero initial
momentum \cite{pytorch2026muon}:
\begin{lstlisting}[style=pythonpaper]
import torch
from torch.utils.data import DataLoader

loader = DataLoader(
    dataset, batch_size=1, shuffle=True, drop_last=False,
)

scale = max(1.0, W.shape[0] / W.shape[1])**0.5
optimizer = torch.optim.Muon(
    [W], lr=eta0 / scale, momentum=beta,
    nesterov=False, weight_decay=0.0, ns_steps=5,
    adjust_lr_fn="original",
)
\end{lstlisting}
Here \texttt{dataset} is a map-style dataset containing the $N$ fixed components. The setting \texttt{batch\_size=1} selects one component per update, as in \cref{alg:rr-muon}. With \texttt{shuffle=True}, each epoch visits every component once in a permuted order. 
The parameter \texttt{beta} is the momentum coefficient $\beta$, and \texttt{eta0} is the initial stepsize $\eta_0$. Dividing \texttt{eta0} by \texttt{scale} compensates for the optimizer's matrix-shape scaling. This call retains the finite-step Newton--Schulz iteration, which \cref{alg:rr-muon} idealizes by the exact polar map.

\section{Experimental Details}\label[appendix]{app:experiments}

This appendix records the settings of \cref{sec:experiments}. In all experiments, the momentum buffer starts from zero and is kept across epochs, and $\beta$ is constant within each run. RR draws a fresh uniformly random permutation in each epoch, SO draws one uniformly random permutation and reuses it, IG uses the order $1,\ldots,N$, and sampling with replacement draws each component independently and uniformly at random. Unless stated otherwise, reported values are medians over random seeds. The initialization is the same across seeds. Each subsection first gives the setting common to all experiments and then the settings of each experiment.

\subsection{MLP on MNIST}\label[appendix]{app:exp-mnist}

\textbf{Common setting.} The $S=60000$ training images are flattened and scaled to $[0,1]$. The network $x\mapsto W_2\tanh(W_1x)$ has $W_1\in\R^{64\times784}$ and
$W_2\in\R^{10\times64}$ without biases, initialized with independent $\mathcal{N}(0,1/784)$ and $\mathcal{N}(0,1/64)$ entries from a fixed seed.
The images are split into $N=S/B$ disjoint batches $\mathcal{B}_1,\ldots,\mathcal{B}_N$ of size $B$, and batch $i$ defines the component
\[
    f_i(W_1,W_2)=\frac{1}{B}\sum_{j\in\mathcal{B}_i}
        \ell_{\mathrm{CE}}\prt{W_2\tanh(W_1x_j),y_j}
        +\frac{\lambda}{2}\prt{\normf{W_1}^2+\normf{W_2}^2},
    \qquad i=1,\ldots,N,\quad \lambda=10^{-4},
\]
where $\ell_{\mathrm{CE}}$ is the softmax cross-entropy. Problem~\eqref{eq:P} is then
\[
    f(W_1,W_2)=\frac{1}{N}\sum_{i=1}^{N}f_i(W_1,W_2)
    =\frac{1}{S}\sum_{j=1}^{S}\ell_{\mathrm{CE}}\prt{W_2\tanh(W_1x_j),y_j}
        +\frac{\lambda}{2}\prt{\normf{W_1}^2+\normf{W_2}^2},
\]
which does not depend on $B$.

Muon applies the exact polar map to each layer's momentum buffer. Unless stated otherwise, each run has $T=2\times10^4$ iterations and the stepsize is set to $\eta=3\times10^{-4}$ for every $\beta$, $B$, and sampling scheme.

\textbf{Sampling schemes.} The sampling-scheme comparison uses $B=100$ ($N=600$), $\beta\in\crki{0.9,0.99}$, and three seeds for RR, SO, and sampling with
replacement. IG is deterministic and run once. The random partition is a uniformly random split. The sorted partition sorts the images by label, with
ties broken at random, before splitting. 

\textbf{Required momentum and the number of components.} The batch sizes are $B\in\crki{25,50,100,200,400,800}$, which gives $N=S/B\in \mathcal{N}_{\text{MNIST}}:= \crki{2400, 1200,600,300,150,75}$. We shuffle the images once and, for each $B$, split them into consecutive batches of size $B$. The grid is given by 
\begin{align*}
    \crk{(N, \beta): N\in\mathcal{N}_{\text{MNIST}},\ \beta\in\crki{1-10^{-k/4}: k=0,1,\ldots,12} }.
\end{align*}
Each point of the grid is run with three random seeds, and the heatmap and the level curves in \Cref{fig:exp-batch} (left) show the median final training loss $f(W_T)$ over the seeds.

\textbf{Tuning suggestion (i).} Both panels of \Cref{fig:exp-tuning-budget} use $B=100$ and IG on the sorted partition. The left figure compares all combinations of $\beta\in\crki{0.9,0.99}$ and $\eta\in\crki{3\times10^{-4},10^{-4},3\times10^{-5}}$, each run for $T=6/\eta$ iterations. For the right panel, we first run the baseline with $\beta=0.9$ and $\eta=3\times10^{-4}$ for $4800$ updates, by which point the training loss has reached a plateau. We then continue from this point with different choices of $\beta$ and $\eta$. Specifically, we either keep $\beta=0.9$ with $\eta\in\crki{3\times10^{-4},10^{-4},3\times10^{-5}}$ or increase $\beta$ to $0.99$ with $\eta\in\crki{3\times10^{-4},10^{-4}}$.

\textbf{Tuning suggestion (ii).} We use $B=100$, RR on the random partition, $T=2\times10^4$, and three seeds. With $\beta_k=1-10^{-k/2}$ and $\eta_l=10^{l/4-5}$, the grid in \Cref{fig:exp-tuning} (left) is
\begin{align*}
    \crk{(\beta_k,\eta_l): k=2,3,\ldots,6,\ l=4,5,\ldots,16}.
\end{align*}

\subsection{ResNet-18 on CIFAR-10}\label[appendix]{app:exp-cifar}

\textbf{Common setting.} The $50000$ training images are
standardized per channel, without data augmentation, and split uniformly at random into $N=400$ disjoint batches $\mathcal{B}_1,\ldots,\mathcal{B}_N$ of size $B=125$. Batch $i$ defines
the component
\[
    f_i(W)=\frac{1}{B}\sum_{j\in\mathcal{B}_i}
        \ell_{\mathrm{CE}}\prt{h_W(x_j;\mathcal{B}_i),y_j}
        +\frac{\lambda}{2}\sum_{l}\normf{W_l}^2,
    \qquad i=1,\ldots,N,\quad \lambda=10^{-4},
\]
where $W=(W_l)_l$ collects the matrices of the network and
$h_W(x_j;\mathcal{B}_i)$ is its output for image $x_j$ with BatchNorm in
training mode. 
%Since BatchNorm normalizes each image with the statistics of its batch, $f = \frac{1}{N}\sum_{i=1}^{N} f_i$ depends on the partition, unlike for the MLP.

We use the CIFAR variant of ResNet-18, with a $3\times3$
stride-one first convolution and no max pooling. To update every
trainable parameter by Muon, we remove the affine parameters of BatchNorm and the bias of the classifier, and reshape each convolution kernel, which has $C_{\mathrm{out}}$ output channels, $C_{\mathrm{in}}$ input channels, and size $k\times k$, into a $C_{\mathrm{out}}\times C_{\mathrm{in}}k^2$ matrix. The initialization uses a fixed seed.

All parameters are updated by Muon with non-Nesterov momentum $\beta=0.95$, the default coefficient in \cite{jordan2024muon}. The polar map is approximated by five Newton--Schulz iterations.

Each matrix $W_l\in\R^{m_l\times n_l}$ is updated with the constant stepsize $0.02\max\crki{1,m_l/n_l}^{1/2}$, where $m_l$ and $n_l$ are its numbers of rows and columns, e.g., $m_l=C_{\mathrm{out}}$ and $n_l=C_{\mathrm{in}}k^2$ for a convolution. Here $0.02$ is the default stepsize and $\max\crki{1,m_l/n_l}^{1/2}$ the shape factor of the reference implementation~\cite{jordan2024muon}. To tune the hyperparameters $\lambda$ and $\beta$, we run RR for ten epochs with the stepsize above and $(\lambda,\beta)\in\crki{(10^{-4},0.95),(5\times10^{-4},0.95),(10^{-4},0.99)}$, and choose $\lambda=10^{-4}$ and $\beta=0.95$, which give the lowest training loss and the highest test accuracy among the pairs.

\textbf{Sampling schemes.} Each run has $30$ epochs ($T=12000$). Every two
epochs, we evaluate the full training loss as the average of all $400$
components in training mode, without updating the BatchNorm running
statistics. RR, SO, and sampling with replacement use three seeds, and IG is run once.
\Cref{fig:exp-cifar-gpt2} (left) shows the median over seeds and the range.

\subsection{GPT-2 on FineWeb-Edu}\label[appendix]{app:exp-gpt2}

\textbf{Common setting.} We stream the \texttt{sample-10BT} subset of FineWeb-Edu~\cite{penedo2024fineweb}, append the end-of-text token, and keep the first
$5\times10^6$ tokens for validation and the next $3\times10^8$ for training.\footnote{The tokenization, the construction of the sequences, the padded vocabulary, the weight tying, and the initialization follow A.~Karpathy, \texttt{build-nanogpt}, files \texttt{fineweb.py} and \texttt{train\_gpt2.py}: \url{https://github.com/karpathy/build-nanogpt}. There, the end-of-text token precedes each document, and $10^8$ tokens are used for validation and about $10^{10}$ for training.}
The training tokens are cut into $292968$ sequences of $1025$ tokens, in which consecutive sequences share one token, so that each sequence gives
$1024$ next-token predictions. We drop the last eight and split the remaining $292960$ sequences once into $N=9155$ disjoint batches
$\mathcal{B}_1,\ldots,\mathcal{B}_N$ of $B=32$ sequences. Batch $i$ defines the component
\[
    f_i(W)=\frac{1}{1024B}\sum_{j\in\mathcal{B}_i}\sum_{k=1}^{1024}
        -\log p_W\prt{s_{j,k+1}\mid s_{j,1},\ldots,s_{j,k}},
    \qquad i=1,\ldots,N,
\]
the mean cross-entropy of the $32\times 1024$ next-token predictions of batch $i$, where $s_j=(s_{j,1},\ldots,s_{j,1025})$ is sequence $j$ and $p_W$ is the distribution predicted by the model. 

We use GPT-2 small~\cite{radford2019language} with the vocabulary padded to $50304$. To update every trainable parameter by Muon, we use separate query, key, and value matrices and remove the biases of the linear layers and the affine parameters of LayerNorm. With the output layer tied to the token embedding, the trainable parameters are $74$ matrices with $124.4$M entries in total, including the token and position embeddings.
All entries are drawn independently from $\mathcal{N}(0,0.02^2)$, except those of the output projections in the attention and MLP blocks, which are drawn from $\mathcal{N}(0,(0.02/\sqrt{24})^2)$. The initialization uses a fixed seed.

All matrices, including the embeddings, are updated by Muon with non-Nesterov momentum $\beta=0.95$ and five Newton--Schulz iterations, computed in bfloat16. The stepsize is constant, $\eta\max\crki{1,m/n}^{1/2}$ for an $m\times n$ parameter with $\eta=0.005$. 
For each candidate stepsize $\eta\in\crki{0.0025,0.005,0.0075,0.01,0.02,0.04}$, we run $1000$ RR updates and select the value with the lowest estimated training loss.

\textbf{Sampling schemes.} Each run has three epochs ($T=27465$). We estimate the training loss as the average of $128$ fixed components drawn uniformly without replacement. RR, SO, and sampling with replacement use three seeds, and IG is run once. 

\textbf{Tuning suggestion (ii).} We run RR for one epoch ($T=N=9155$) with one seed and estimate the final training loss on $512$ instead of $128$ components. With $\beta_k=1-10^{-k/2}$ and $\eta_l=5\times10^{-4}\cdot2^{l}$, the grid in \Cref{fig:exp-tuning} (right) is
\begin{align*}
    \crk{(\beta_k,\eta_l): k=2,3,\ldots,6,\ l=0,1,\ldots,6}.
\end{align*}

\subsection{The smooth convex counterexample}\label[appendix]{app:exp-counterexample}

\textbf{Common setting.} We use the components from \eqref{eq:huber-matrix} with $N=4K+1\in\crki{5,9,17,33,65,129,257,513}$, i.e., $K=2^q$ for $q=0,\ldots,7$. By \cref{app:smooth-counterexample-proof}, the iterates of Muon from $W_0=M_{-1}=0$ are $W_t=x_tE_{11}$, so $f(W_t)-f^*=g(x_t)-g(x_*)$ and $\norm{\nabla f(W_t)}_*=|g'(x_t)|$, where $x_*=-1-1/N$ is the minimizer of $g$. Hence, we run Muon under RR, SO, and IG through the scalar recursion \eqref{eq:scalar-muon}, in which the exact polar map reduces to the sign, from $x_0=m_{-1}=0$ with the stepsizes in \eqref{eq:huber-steps} for $p=0.6$. 
Since $f_2=\cdots=f_N$, a permutation affects the run only through the position of $f_1$ in each epoch: RR draws this position uniformly in every
epoch, SO draws it once, and IG places $f_1$ first. 
% We also run the
% $2\times2$ embedding \eqref{eq:huber-matrix} with an SVD-based polar map for
% $N=17$, $\beta=0.95$, and $300$ epochs of RR. It reproduces the scalar
% iterates exactly, and the entries other than $(1,1)$ remain zero.

\begin{figure}[!htbp]
    \centering
    \begin{minipage}[b]{0.49\linewidth}
        \centering
        \includegraphics[width=\linewidth]{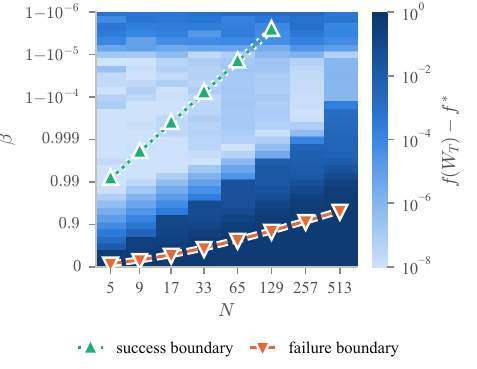}
    \end{minipage}\hfill
    \begin{minipage}[b]{0.49\linewidth}
        \centering
        \includegraphics[width=\linewidth]{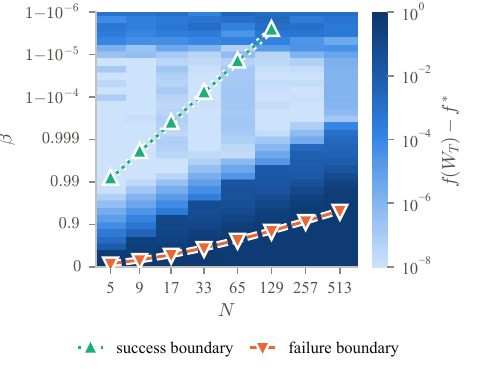}
    \end{minipage}
    \caption{Final gap $f(W_T) - f^*$ on the counterexample under SO (left) and IG (right).}
    \label{fig:exp-counterexample-orders}
\end{figure}

\textbf{Required momentum and the number of components.} The grid is given by
\begin{align*}
    \crk{(N,\beta): N\in\crki{5,9,17,\ldots,513}, \beta\in\crk{1-10^{-k/6}: k=0,1,\ldots,36}}.
\end{align*}
Each point of the grid is run for $T=10^6$ iterations, regardless of $N$. RR and SO use five seeds, and IG is run once. The heatmap in \Cref{fig:exp-counterexample} (left) shows the median final gap $f(W_T)-f^*$ over the seeds under RR, and \Cref{fig:exp-counterexample-orders} shows the same for SO and IG. The theoretical failure boundary is $(4K+2)\beta^K=1$, below which \cref{prop:smooth-n1} applies. The theoretical success boundary is $2N(N+1)\sqrt{N-1}(1-\beta)=1$, above which the residual $R$ in \eqref{eq:residual} is below the floor $1/N$ of \eqref{eq:huber-failure}. 

The trajectories in \Cref{fig:exp-counterexample} (right) use RR with $N\in\crki{9,33,129,513}$, $\beta\in\crki{0.9,0.99,0.999,0.9999}$, $T=10^6$, and ten seeds. We record $f(W_t)-f^*$ and show the median over the seeds and the band between the $10\%$ and $90\%$ quantiles.

\FloatBarrier

\bibliographystyle{siamplainnat}
\bibliography{reference}

\end{document}